\documentclass[11pt, eqno]{article}
\usepackage{amssymb,amsmath,latexsym}
\usepackage{epsfig}
\usepackage{color}

\newtheorem{thm}{Theorem}[section]
\newtheorem{cor}[thm]{Corollary}
\newtheorem{lemma}[thm]{Lemma}
\newtheorem{prop}[thm]{Proposition}
\newtheorem{defn}[thm]{Definition}

\newtheorem{remark}[thm]{Remark}
\numberwithin{equation}{section}

\def\pf{{\medskip\noindent {\bf Proof. }}}
\def\pff{{\medskip\noindent {\bf Proof }}}
\def\qed{{\hfill $\Box$ \bigskip}}

\def\sA {{\cal A}} \def\sB {{\cal B}} 
 \def\sE {{\cal E}} \def\sF {{\cal F}}

\def\sM {{\cal M}}

  \def\bR {{\mathbb R}}

\def\R {{\mathbb R}} \def\RR {{\mathbb R}}

\def\wt{\widetilde}
\def\wh{\widehat}

\def\E{{\mathbb E}}
\def\P{{\mathbb P}}

\def\bea{\begin{align*}}
\def\eea{\end{align*}}
\def\bee{\begin{equation}}
\def\eee{\end{equation}}

\def\beq{\begin{equation}}
\def\eeq{\end{equation}}

\def\eps{\varepsilon}

\def\wh{\widehat}

\def\1{{\bf 1}}

\makeatletter
\@addtoreset{equation}{section}

\makeatother

\begin{document}
\bibliographystyle{plain}

\title{\Large \bf
Sharp Green
Function Estimates for $\Delta + \Delta^{\alpha/2}$ in $C^{1,1}$ Open Sets
and Their Applications
}

\author{{\bf Zhen-Qing Chen}\thanks{Research partially supported
by NSF Grant  DMS-0906743.}, \quad {\bf Panki Kim}\thanks{This work
was supported by Basic Science Research Program through the National
Research Foundation of Korea(NRF) grant funded by the Korea
government(MEST)(2009-0093131).}, \quad {\bf Renming Song} \quad and
\quad {\bf Zoran Vondra\v{c}ek}\thanks{Research partially supported
by the MZOS grant 037-0372790-2801 of the Republic of Croatia.}}

 \date{(December 11, 2009)}

\maketitle

\begin{abstract}
We consider a family of pseudo differential operators $\{\Delta+
a^\alpha \Delta^{\alpha/2}; \ a\in [0, 1]\}$ on $\R^d$ that evolves
continuously from $\Delta$ to $\Delta + \Delta^{\alpha/2}$, where
$d\geq 1$ and $\alpha \in (0, 2)$. It gives rise to a family of
L\'evy processes \{$X^a, a\in [0, 1]\}$, where $X^a$ is the sum of a
Brownian motion and an independent symmetric $\alpha$-stable process
with weight $a$. Using a recently obtained uniform boundary Harnack
principle with explicit decay rate, we  establish sharp bounds for
the Green function of the process $X^a$ killed upon exiting a
bounded $C^{1,1}$ open set $D\subset\R^d$. As a
 consequence, we identify the Martin boundary of $D$ with respect to
$X^a$ with its Euclidean boundary. Finally, sharp Green function
estimates are derived for certain L\'evy processes which can be
obtained as perturbations of $X^a$.
\end{abstract}

\bigskip
\noindent {\bf AMS 2000 Mathematics Subject Classification}:
Primary 31A20,  31B25, 60J45; Secondary   47G20,  60J75, 31B05

\bigskip\noindent
{\bf Keywords and phrases}: Green function estimates, boundary
Harnack principle, harmonic functions, fractional Laplacian,
Laplacian, symmetric $\alpha$-stable process, Brownian motion,
perturbation

\bigskip
\section{Introduction}
Discontinuous Markov processes have been intensively studied in
recent years due to their importance both in theory and
applications. In contrast to the diffusion case, the infinitesimal
generator of a discontinuous Markov process in $\R^d$ is a non-local
(or integro-differential) operator. Most of the recent studies have
concentrated on discontinuous Markov processes (and corresponding
integro-differential operators) that do not have a diffusion
component. See
\cite{BBKRSV, C0}
and the references therein for a
summary of some of these recent results from the probability
literature. For recent progress in PDE literature, we refer the
readers to \cite{CSS, CaS, CV}.

However, Markov processes with both  diffusion and jump components
are needed in many situations, like in finance and control
 theory. See,  for example, \cite{JKC,
 MP, OS}. On the other hand,
the fact that such a process $X$ has both diffusion and jump
components is also the source of many technical difficulties in
investigating the potential theory of $X$. The main difficulty in
studying $X$ stems from the fact that it runs on two different
scales: on the small scale the diffusion part dominates, while on
the large scale the jumps take over. Another difficulty is
encountered at the exit of $X$ from an open set: for diffusions, the
exit is through the boundary, while for the pure jump processes,
typically the exit happens by jumping out from the open set. For the
process $X$, both cases will occur which makes the process $X$ much
more difficult to study.

Despite these difficulties, significant progress has been made in
the last few years in understanding the potential theory of
discontinuous Markov processes with both diffusion and jump
components. Green function estimates (for the whole space) and the
Harnack inequality for some  processes with both diffusion and jump
components were established in \cite{RSV, SV05}. The parabolic
Harnack inequality and heat kernel estimates were studied in
\cite{SV07} for the independent sum of a Brownian motion and a
symmetric stable process, and in \cite{CK08} for much more general
diffusions with jumps. Moreover, an a priori H\"older estimate is
established in \cite{CK08} for bounded parabolic functions. Very
recently, the boundary Harnack principle for some one-dimensional
L\'evy processes with both diffusion  and jump components was
studied in \cite{KSV09}, where sharp estimates on Green functions of
bounded open sets of $\R$ were also established.
Most recently, a boundary Harnack principle with explicit decay rate
for nonnegative harmonic functions of the independent sum of a
Brownian motion and a symmetric stable process in $C^{1,1}$ open
sets in $\R^d$ was obtained in \cite{CKSV}.

The main goal of this paper is to use the boundary Harnack principle
obtained in \cite{CKSV} to establish sharp Green function estimates
in $C^{1,1}$ open sets for the L\'evy processes that are
independent sums of Brownian motions and symmetric stable processes.

Let us now fix the notation and state the main result of this paper.
Throughout this paper, we assume that $d \ge 1$ is an integer and
$\alpha\in (0, 2)$. Let $X^0=(X^0_t,\, t\ge 0)$ be a Brownian motion
in $\R^d$ with generator $\Delta =\sum_{i=1}^d
\frac{\partial^2}{\partial x_i^2}$, and let $Y=(Y_t,\, t\ge 0)$ be
an independent (rotationally) symmetric $\alpha$-stable process in
$\R^d$. For $a>0$ we define the process $X^a=(X^a_t,\, t\ge 0)$ by
$X^a_t=X^0_t+a Y_t$, called \emph{the independent sum of a Brownian
motion and a symmetric $\alpha$-stable process with weight $a$}.
This process, although very specific, serves as a test case for
general Markov processes with both diffusion and jump components.

Let $D$ be a $C^{1,1}$ open set in $\R^d$, let $X^{a,D}$ be the
process $X^a$ killed upon exiting $D$ and let $G_D^a(x,y)$ denote
the Green function of $X^{a,D}$ (for precise definitions see Section
2). Our main goal is to establish sharp two-sided estimates for
$G_D^a(x,y)$. Let $\delta_D(x)$ denote the Euclidean distance
between the point $x\in D$ and the boundary $\partial D$. The main
result of this paper is the following theorem.
Here and in the sequel,
for $a, b\in \bR$, $a\wedge b:=\min \{a,
b\}$ and $a\vee b:=\max\{a, b\}$
Define for $d\geq 3$
and $a>0$,
$$
g_D^a (x, y) := \begin{cases} \frac{1} {|x-y|^{d-2}} \left(1\wedge
\frac{  \delta_D(x) \delta_D(y)}{ |x-y|^{2}}\right)
\quad &\hbox{when } x, y \hbox{ are in the same component of } D, \\
\frac{a^\alpha} {|x-y|^{d-2}} \left(1\wedge \frac{  \delta_D(x)
\delta_D(y)}{ |x-y|^{2}}\right) \quad &\hbox{when } x, y \hbox{ are
in different components of } D;
\end{cases}
$$
for $d=2$ and $a>0$,
\begin{equation}\label{e:1.2}
g_D^a (x, y): = \begin{cases} \log\left(1+\frac{  \delta_D(x)
\delta_D(y)}{ |x-y|^{2}}\right)
\quad &\hbox{when } x, y \hbox{ are in the same component of } D, \\
a^\alpha \log\left(1+\frac{  \delta_D(x) \delta_D(y)}{
|x-y|^{2}}\right) \quad &\hbox{when } x, y \hbox{ are in different
components of } D;
\end{cases}
\end{equation}
and for $d=1$ and $a>0$,
\begin{equation}\label{e:1.21d}
g_D^a (x, y): = \begin{cases} \left(\delta_D(x)
\delta_D(y)\right)^{1/2}\wedge\frac{ \delta_D(x) \delta_D(y)}{|x-y|}
\quad &\hbox{when } x, y \hbox{ are in the same component of } D, \\
a^\alpha \big(\left(\delta_D(x) \delta_D(y)\right)^{1/2}\wedge\frac{
\delta_D(x) \delta_D(y)}{|x-y|}\big)\quad &\hbox{when } x, y \hbox{
are in different components of } D.
\end{cases}
\end{equation}

\begin{thm}\label{t-main-green}
Let $M>0$. Suppose that $D$ is a bounded $C^{1,1}$
 open set in $\R^d$. There exists
 $C=C(D, M, \alpha)>1$ such that for
all $x, y \in D$ and all $a \in (0,M]$
\begin{equation}\label{e:1.3}
C^{-1} \, g_D^a(x, y) \leq G_D^a(x, y) \leq C\, g_D^a(x, y) .
\end{equation}
\end{thm}

Note that the above estimates are uniform in $a\in (0, M]$. In case
$d=1$, a (non-uniform) estimate is covered by \cite{KSV09}. Letting
$a\downarrow 0$  in \eqref{e:1.3} recovers the Green function
estimates for Brownian motion killed upon exiting $D$; for the
latter, see \cite[p.~182]{CZ} for $d=2$ and \cite{Zh} for $d\geq 3$,
respectively.

The rest of the paper is organized as follows. Section 2 gives
preliminary and background materials.  Theorem \ref{t-main-green} is
proved in Sections 3, 4 and 5. The proof of the theorem in the case
$d\geq 3$ is by now quite standard. Once the interior estimates are
established, the full estimates in connected $C^{1,1}$ open sets
follow from the boundary Harnack principle by the method developed
by Bogdan \cite{Bo1} and Hansen \cite{H}. However this method is not
applicable when $d\le 2$ since Brownian motion is recurrent in this
case. When $d=2$, we use a capacitary argument and some recent
results on subordinate killed Brownian motions, which are given in
Section 4. The case $d=1$ is dealt with in Section 5, where we
follow the arguments of \cite{KSV09}. In Section 6, using the
boundary Harnack principle and our Green function estimates, we show
that both the Martin and the minimal Martin boundary of the process
$X^{a,D}$ can be identified with the Euclidean boundary when $D$
 is a bounded $C^{1,1}$ open set. In the last section, we extend our
results on $X^a$  to symmetric L\'evy processes that can be obtained
from $X^a$ through certain perturbations. In particular, for every
$m>0$, we obtain sharp Green function estimates of $\Delta + m -
(m^{2/\alpha} -\Delta)^{\alpha}$ in any bounded $C^{1,1}$ open set
with zero exterior condition. The process corresponding to $\Delta +
m - (m^{2/\alpha} -\Delta)^{\alpha}$ is a L\'evy process that is the
independent sum of a Brownian motion and a relativistic
$\alpha$-stable process with mass $m$.

\medskip

Throughout this paper, we use the capital letters $C_1,C_2, \cdots $
to denote constants in the statement of results, and their labeling
will be fixed. The lowercase constants $c_1, c_2, \cdots$ will
denote generic constants used in proofs, whose exact values are not
important and can change from one appearance to another. The
labeling of the constants $c_1, c_2, \cdots$ starts anew in every
proof. The dependence of the constant $c$ on the dimension $d$ and
$\alpha \in (0, 2)$ may not be mentioned explicitly. The constant
$\alpha \in (0, 2)$ will be fixed throughout this paper.
We will use ``$:=$" to denote a
definition, which is read as ``is defined to be". $B(x,r)$
denotes the open ball in $\R^d$ centered at $x$ with radius $r>0$.
Recall that for any $x\in D$, $\delta_D(x)$ denotes the distance
between $x$ and $\partial D$, and for $a, b\in \bR$,
$a\wedge b:=\min \{a, b\}$ and $a\vee b:=\max\{a, b\}$.
   We will use $\partial$ to
denote a cemetery point and for every function $f$, we extend its
definition to $\partial$ by setting $f(\partial )=0$. Lebesgue
measure in $\bR^d$ will be denoted by $dx$. For a Borel set
$A\subset \bR^d$, we also use $|A|$ to denote its Lebesgue measure.

\section{Preliminaries}

A (rotationally) symmetric $\alpha$-stable process $Y=(Y_t, t\geq 0,
\P_x, x\in \bR^d)$ in $\bR^d$ is a L\'evy process with the
characteristic exponent $|\xi|^{\alpha}$, i.e.,
$$
\E_x \left[ e^{i\xi\cdot(Y_t-Y_0)} \right]\,=\,e^{-t|\xi|^{\alpha}}
\qquad \hbox{for every } x\in \bR^d \hbox{ and }  \xi\in \bR^d.
$$
The infinitesimal generator of $Y$
 is the fractional Laplacian $\Delta^{\alpha /2}$, which is a
prototype of non-local operators. The fractional Laplacian
can be written in the form
$$
\Delta^{\alpha /2} u(x)\, =\, \lim_{\eps \downarrow 0}\int_{\{y\in
\bR^d: \, |y-x|>\eps\}} (u(y)-u(x)) \frac{\sA (d, -\alpha)
}{|x-y|^{d+\alpha}}\, dy,
$$
where $ {\cal A}(d, -\alpha):= \alpha2^{\alpha-1}\pi^{-d/2}
\Gamma(\frac{d+\alpha}2) \Gamma(1-\frac{\alpha}2)^{-1}. $ Here
$\Gamma$ is the Gamma function defined by $\Gamma(\lambda):=
\int^{\infty}_0 t^{\lambda-1} e^{-t}dt$ for every $\lambda > 0$.

Suppose $X^0$ is a Brownian motion in $\R^d$ with generator
$\Delta=\sum_{i=1}^d \frac{\partial^2}{\partial x_i^2}$, and $Y$ is
a symmetric $\alpha$-stable process in $\R^d$. Assume that $X^0$ and
$Y$ are independent. For any $a>0$, we define the process
$X^a=(X^a_t, \, t\ge 0)$ by $X_t^a:=X^0_t+ a Y_t$. As already
mentioned, the process $X^a$ is called the independent sum of the
Brownian motion $X^0$ and the symmetric $\alpha$-stable process  $Y$
with weight $a$. It is a L\'evy process with the characteristic
exponent $\Phi^a(\xi)=|\xi|^2+a^{\alpha}|\xi|^{\alpha}$, $\xi\in
\bR^d$, and its infinitesimal generator is $\Delta+a^\alpha
\Delta^{\alpha/2}$. The process $X^a$ has a jointly continuous
transition density that will be denoted by $p^a(t,x,y)$. From the
Chung-Fuchs criterion (see  \cite[Theorem I.17]{B}), it easily
follows that, when $a>0$, $X^a$ is transient if and only if
$\alpha<d$, while it is well known that $X^0$ is transient if and
only if $d\ge 3$.

There is another representation of the process $X^a$ which will be
useful in Sections 3, 4 and 5. It can be obtained by subordinating
$X^0$ with an independent subordinator $T^a_t:=t+a^{2} T_t$ where
$T=(T_t,\, t\ge 0)$ is an $\alpha/2$-stable subordinator, i.e., the
processes $(X^a_t, \,t \ge 0)$ and $(X^0_{T^a_t}, \,t \ge 0)$ have
the same distribution. Note that the Laplace exponent
of
 $T^a$ is $\phi^a(\lambda)=\lambda +a^{\alpha}\lambda^{\alpha/2}$.
Let $\sM_{\alpha/2}(t):=\sum_{n=0}^{\infty}(-1)^n
t^{n\alpha/2}/\Gamma(1+n\alpha/2)$. It follows by a straightforward
integration that
$$
\int_0^{\infty}e^{-\lambda t}\sM_{1-\alpha/2}(a^{2\alpha/(2-\alpha)}
t)\, dt=\frac{1}{\phi^a(\lambda)}\, ,
$$
which shows that the potential density $u^a$ of the subordinator $T^a$
is given by
\begin{equation}\label{e:pd4sub}
u^a(t)=\sM_{1-\alpha/2}(a^{2\alpha/(2-\alpha)} t)\, .
\end{equation}
Since, for any $a>0$, $\phi^a$ is a complete Bernstein function, we
know that (see, for instance, \cite{RSV}) $u^a(\cdot)$ is a
completely monotone function. In particular, $u^a(\cdot)$ is a
decreasing function. Since $u^a(t)=u^1(a^{2\alpha/(2-\alpha)} t)$,
we see that $a\mapsto u^a(t)$ is a decreasing function. Moreover,
since the drift of $T^a$ is equal to 1, we have that $u^a(0+)=1$ and
so
\begin{equation}\label{e:le}
u^a(t) \leq
1 \quad \text{ for }t>0.
\end{equation}

The L\'evy measure of $X^a$ has a density with respect to the
Lebesgue measure given by
\begin{equation}\label{e:j}
J^a (x, y):=j^a(y-x):=j^a (|y-x|)=a^\alpha \sA(d,
-\alpha)|x-y|^{-(d+\alpha)},
\end{equation}
which is called the L\'evy intensity of $X^a$. It determines a
L\'evy system for $X^a$, which describes the jumps of the process
$X^a$: For any non-negative measurable function $f$ on $\bR_+ \times
\bR^d\times \bR^d$ with $f(s, x, x)=0$ for all $s>0$ and  $x\in
\bR^d$, and stopping time $T$ (with respect to the filtration of
$X^a$),
\begin{equation}\label{e:levy}
\E_x \left[\sum_{s\le T} f(s,X^a_{s-}, X^a_s) \right]= \E_x \left[
\int_0^T \left( \int_{\bR^d} f(s,X^a_s, y) J^a(X^a_s,y) dy \right)
ds \right]
\end{equation}
(see, for example, \cite[Proof of Lemma 4.7]{CK} and \cite[Appendix
A]{CK2}).

The quadratic form $(\sE^a, \sF)$ associated with the generator
$\Delta + a^\alpha \Delta^{\alpha/2}$ of $X^a$ is given by
$$
\sF=W^{1,2}(\R^d):=\left\{u\in L^2(\R^d; dx): \,
 \frac{\partial u}{\partial x_i} \in L^2(\R^d; dx)
 \, \hbox{ for every } \,
  1\le i \le d \right\}
$$
and for  $u, v\in \sF$,
$$
\sE^a(u, v) = \int_{\R^d} \nabla u(x) \cdot \nabla v(x) \, dx +
\frac{1}{2} \int_{\R^d\times \R^d} (u(x)-u(y))(v(x)-v(y)) \frac{
\sA(d, -\alpha)\, a^\alpha }{|x-y|^{d+\alpha}} dxdy .
$$
In probability theory, the quadratic form $(\sE^a, W^{1,2}(\R^d))$
is called the Dirichlet form of $X^a$.
Let
$\sE^a_1(u,
u):=\sE^a(u, u)+\int_{\R^d} u(x)^2 dx$. Note that for every $a>0$,
there is a positive constant $c=c(a, d, \alpha)\geq 1$  so that
$$
\int_{\R^d}  \left( |\nabla u(x)|^2 + u(x)^2 \right) dx \leq
\sE^a_1(u, u) \leq c \, \int_{\R^d} \left( |\nabla u(x)|^2 + u(x)^2
\right) dx \qquad \hbox{for } u\in
W^{1,2}(\R^d).
$$
Thus the processes $X^a$, $a\ge0$, share the same family of sets
having zero capacity.

For any open set $D\subset \bR^d$, $\tau^a_D:=\inf\{t>0: \,
X^a_t\notin D\}$ denotes the first exit time from $D$ by $X^a$. We
denote by  $X^{a,D}$ the subprocess of $X^a$ killed upon leaving
$D$. The infinitesimal generator of $X^{a,D}$ is $(\Delta+ a^\alpha
\Delta^{\alpha/2})|_D$. It is known (see \cite{CK08}) that $X^{a,D}$
has a continuous transition density $p^a_D(t, x, y)$ with respect to
the Lebesgue measure.

\begin{defn}\label{D:1.1}  \rm A real-valued function $u$ defined on
$\R^d$ is said to be
\begin{description}
\item{(1)}
harmonic in $D\subset \R^d$ with respect to $X^a$ if for every open
set $B$ whose closure is a compact subset of $D$,
\begin{equation}\label{e:har}
\E_x \left[ \big| u(X^a_{\tau^a_{B}})\big| \right]<\infty \quad
\hbox{and} \quad u(x)= \E_x \left[ u(X^a_{\tau^a_{B}})\right] \qquad
\hbox{for every } x\in B;
\end{equation}

\item{(2})
regular harmonic in $D\subset \R^d$ with respect to $X^a$ if it
is harmonic in $D$ with respect to $X^a$ and
$$
u(x)= \E_x\left[u(X^a_{\tau^a_{D}})\right]\quad  \hbox{for every } x\in D;
$$

\item{(3}) harmonic for $X^{a,D}$ if it is harmonic for $X^a$ in $D$ and
vanishes outside $D$;

\item{(4)}
superharmonic in $D\subset \R^d$ with respect to $X^a$ if for every open
set $B$ whose closure is a compact subset of $D$,
\begin{equation}\label{e:suphar}
\E_x \left[ \big| u(X^a_{\tau^a_{B}})\big| \right]<\infty \quad
\hbox{and} \quad u(x)\ge  \E_x \left[ u(X^a_{\tau^a_{B}})\right] \qquad
\hbox{for every } x\in B.
\end{equation}
\end{description}
\end{defn}

\medskip

It follows from \cite{CK08} that every harmonic function $u$ in $D$
with respect to $X^a$ is continuous in $D$ and $\int_{\R^d}
|u(y)|(1\wedge |y|^{-(d+\alpha)})dy <\infty$.

Using  the parabolic Harnack inequality from \cite[Theorem
6.7]{CK08} and a scaling argument, the following uniform Harnack
principle was established in \cite{CKSV}.

\begin{prop}[Uniform Harnack principle]\label{uhp}
Suppose that $M>0$. There exists a constant $C_1=C_1(\alpha, M)>0$
such that for any $r\in (0, 1]$, $a\in [0, M]$, $x_0\in \R^d$ and
any function $u$ which is nonnegative in $\R^d$ and harmonic in
$B(x_0, r)$ with respect to $X^a$ we have
$$
u(x)\le C_1 u(y) \qquad \mbox{ for all } x, y\in B(x_0, r/2).
$$
\end{prop}

We recall that an open set $D$ in $\bR^d$ with $d\geq 2$ is said to
be $C^{1,1}$ if there exist a localization radius $R>0$ and a
constant $\Lambda >0$ such that for every $Q\in \partial D$, there
exist a $C^{1,1}$-function $\phi=\phi_Q: \bR^{d-1}\to \bR$
satisfying $\phi(0)=0$, $ \nabla\phi (0)=(0, \dots, 0)$, $\| \nabla
\phi  \|_\infty \leq \Lambda$, $| \nabla \phi (x)-\nabla \phi (y)|
\leq \Lambda |x-y|$, and an orthonormal coordinate system $CS_Q$:
$y=(y_1, \dots, y_{d-1}, y_d)=:(\wt y, \, y_d)$ with its origin at
$Q$ such that
$$
B(Q, R)\cap D=\{ y=(\wt y, y_d)\in B(0, R) \mbox{ in } CS_Q: y_d >
\phi (\wt y) \}.
$$
The pair $(R, \Lambda)$ is called the characteristics of the
$C^{1,1}$ open set $D$. Note that a $C^{1,1}$ open set may be
disconnected. Observe that the distance between any two distinct
connected open components of $D$ is at least $R$. By a $C^{1,1}$
open set in $\bR$ we mean an open set which can be written as the
union of disjoint intervals so that the minimum of the lengths of
all these intervals is positive and the minimum of the distances
between these intervals is positive. Note that a $C^{1,1}$ open set
may be
unbounded.  It is well known that any $C^{1, 1}$ open set $D$
satisfies the uniform exterior ball condition:
There exists $\wt R>0$ such that for every $z\in \partial D$, there
is a ball $B^z$ of radius $\wt R$ such that $B^z\subset (\overline
D)^c$  and $\partial B^z \cap \partial D=\{z\}$. Without loss of
generality, throughout this paper, we assume that the
characteristics $(R, \Lambda)$ of a $C^{1, 1}$ open set satisfies
$R=\wt R$.

Observe that for any $C^{1, 1}$ open set with $C^{1,1}$ characteristics
$(R, \Lambda)$, there exists a  constant $\kappa\in (0,1/2]$,
which depends only on
$(R, \Lambda)$, such that for each $Q \in \partial D$
and $r \in (0, R)$, $D \cap B(Q,r)$
contains a ball $B(A_r(Q),\kappa r)$
of radius $\kappa r$.
 In the rest of paper, whenever we deal with
$C^{1,1}$ open sets, the constants $\Lambda$, $R$ and $\kappa$ will
have the meaning described above.

Let $Q\in \partial D$. We will say that a function $u:\bR^d\to \R$
vanishes continuously on $ D^c \cap B(Q, r)$ if $u=0$ on $ D^c \cap
B(Q, r)$ and $u$ is continuous at every point of $\partial D\cap
B(Q,r)$.

The following theorem is the main result of \cite{CKSV}.

\begin{thm}[Uniform Boundary Harnack Principle]\label{t:main}
Suppose that $M>0$. For any $C^{1, 1}$ open set $D$ in $\bR^d$ with
the characteristics $(R, \Lambda)$, there exists a positive constant
$C_2=C_2(\alpha, d, \Lambda, R, M)$ such that for all $a \in [0,
M]$, $r \in (0, R]$, $Q\in \partial D$ and any nonnegative function
$u$ in $\R^d$ that is harmonic in $D \cap B(Q, r)$ with respect to
$X^{a}$ and vanishes continuously on $ D^c \cap B(Q, r)$, we have
\begin{equation}\label{e:bhp_m}
\frac{u(x)}{u(y)}\,\le C_2\,\frac{\delta_D(x)}{\delta_D(y)}  \qquad
\hbox{for every } x, y\in  D \cap B(Q, r/2).
\end{equation}
\end{thm}

A subset $D$ of $\R^d$ is said to be Greenian for $X^a$ if $X^{a,
D}$ is transient. A Greenian set for $X^0$ will be simply called
Greenian. As mentioned in the second paragraph of Section 2, when
$d\ge 2$ and $a>0$, any non-empty open set $D\subset\bR^d$ is
Greenian for $X^a$; and any non-empty open set in $\bR^d$ is
Greenian when $d\geq 3$. An open set $D\subset \bR^2$ is Greenian if
and only if $D^c$
 is non-polar (or equivalently, has positive capacity).
 In particular, every bounded open
set in $\bR^2$ is Greenian.

For any $a>0$ and any Greenian open subset $D$ of $\R^d$ for $X^a$,
we use $G^a_D(x,y)$ to denote the Green function of $X^{a,D}$, i.e.,
\begin{equation}\label{e:pg}
G^a_D(x, y):=\int_0^\infty p^a_D(t, x, y)dt
\end{equation}
where $p^a_D(t, x, y)$ is the  continuous transition density of
$X^{a,D}$ with respect to the Lebesgue measure. The function
$G_D^a(\cdot, \cdot)$ is finite off the diagonal. It follows
immediately from \eqref{e:pg} that $G^a_D(x, y)$ is a positive
continuous symmetric function off the diagonal of $D\times D$ such
that for any Borel measurable function $f\geq 0$,
$$
\E_x \left[ \int_0^{\tau^a_D} f(X^a_s) ds \right] =\int_D G^a_D (x,
y) f(y) \, dy.
$$
We set $G^a_D$ equal to zero outside $D\times D$. The function
$G^a_D(x, y)$ is also called the Green function of $X^a$ in $D$. For
any $x \in D$, $G^a_D( \, \cdot \, ,x)$ is superharmonic in $D$ with
respect to $X^a$, harmonic in $D\setminus \{x\}$ with respect to
$X^a$ and regular harmonic in $D\setminus \overline{B(x, \eps)}$
with respect to $X^a$ for every $\eps >0$.

Recall that a point $z$ on the boundary $\partial D$ of an open set
$D$ is said to be a regular boundary point for $X^a$ if
$\P_z(\tau^a_D=0)=1$. An open set $D$ is said to be regular for
$X^a$ if every point in $\partial D$ is a regular boundary point for
$X^a$. It is easy to check that every $C^{1,1}$ open set $D$ is
regular for $X^a$ for all $a>0$ and using the argument in the last
paragraph of the proof of
\cite[Theorem 2.4]{CZ}, we conclude that
for any bounded $C^{1,1}$ open set $D$, $G^a_D( \, \cdot \, ,z)$
vanishes continuously on $\partial D$ for every $z \in D$.

Now, as a corollary of the uniform boundary Harnack principle and
the fact that, for any bounded $C^{1,1}$ open set $D$, $G^a_D( \,
\cdot \, ,z)$ vanishes continuously on $\partial D$ for every $z \in
D$, we have the following proposition.

\medskip

\begin{prop}\label{l:Green_L}
Suppose that $M>0$. For any bounded $C^{1, 1}$ open set $D$ in
$\bR^d$ with the characteristics $(R, \Lambda)$, there exists a
positive constant $C_3=C_3(\alpha, d, \Lambda, R, M)>1$ such that
for all $Q \in \partial D$, $r\in (0,  R)$ and  $a \in (0, M]$ we
have
\begin{equation}\label{e:CG_1}
\frac{ G^a_{D}(x,z_1)}{ G^a_{D} (y,z_1)} \, \le \,C_3\, \frac{ G^a_{D}
(x,z_2) }{ G^a_{D} (y,z_2)},
\end{equation}
when $x, y \in D \setminus \overline{B(Q, r)}$ and $z_1, z_2 \in D
\cap B(Q, r/2)$.
\end{prop}

The following scaling property will be used below: If $(X^{a,
D}_t,t\geq 0)$ is the  subprocess in $D$ of the independent sum of a
Brownian motion and a symmetric $\alpha$-stable process in $\R^d$
with weight $a$, then $(\lambda X^{a, D}_{\lambda^{-2} t}, t\geq 0)$
is the subprocess in $\lambda D$ of the independent sum of a
Brownian motion and a symmetric $\alpha$-stable process in $\R^d$
with  weight $a \lambda^{(\alpha-2)/\alpha}$. So for any
$\lambda>0$, we have
\begin{equation}\label{e:scaling}
p^{a\lambda^{(\alpha-2)/\alpha}}_{\lambda D} ( t,  x, y)=
\lambda^{-d} p^{a}_D (\lambda^{-2}t, \lambda^{-1} x, \lambda^{-1} y)
\qquad \hbox{for } t>0 \hbox{ and } x, y \in \lambda D.
\end{equation}
By integrating the above equation with respect to $t$, we get
that when $D$ is Greenian for $X^a$,
\begin{equation}\label{e:scaling1}
G^{a}_D ( x, y) = \lambda^{d-2} G^{a\lambda^{(
\alpha-2)/\alpha}}_{\lambda D} (\lambda x, \lambda y) \qquad
\hbox{for } x, y \in D.
\end{equation}
In particular, for $d=2$, we have
\begin{equation}\label{e:scaling2}
G^{a}_D ( x, y) = G^{a\lambda^{(\alpha-2)/\alpha}}_{\lambda D}
(\lambda x, \lambda y) \qquad \hbox{for } x, y \in D.
\end{equation}

\section{Higher Dimensional Case: $d\ge 3$}

In this section we assume that $d\ge 3$. We will use
$G^a(x,y)=G^a(y-x)
=G^a_{\R^d}(x,y)$ to denote the Green function of $X^a$ in
$\bR^d$.

Recall that $u^a$ is the potential density of the subordinator
$T^a_t=t+a^{2} T_t$ given in \eqref{e:pd4sub}. The Green function
$G^a$ of $X^a$ is also given by the following formula (\cite{RSV})
\beq\label{e:sg} G^a(x)=\int^{\infty}_0(4\pi
t)^{-d/2}e^{-|x|^2/(4t)}u^a(t)dt, \qquad x\in \R^d. \eeq Using this
formula, we can easily see that $G^a$ is radially decreasing and
continuous in $\R^d\setminus \{0\}$.

\begin{lemma}\label{t:Gorigin}
Suppose that $M>0$. For all $a \in [0,M]$, we have
$$
G^M (x) \le G^a(x)\le \frac1{|x|^{d-2}}, \qquad \text{ for all } x
\in \R^d.
$$
\end{lemma}

\pf We have seen that for all $t>0$, $u^a(t) \leq 1$, and the
function $a\mapsto u^a(t)$ is decreasing  on $\R_+$,
cf.~\eqref{e:le} and the text preceding it. The desired inequalities
follow immediately from these properties and \eqref{e:sg}. \qed

\begin{lemma}\label{G_G}
Suppose that $M>0$. There exist $R_1 >0$ and $L_1 \ge 2$ such that
for all $a \in [0,M]$
 \beq \label{e:G_G}
G^a(x) \, \ge \, 2  \, G^a(L_1 x), \qquad \text{for all } |x| < R_1.
 \eeq
\end{lemma}

\pf By \cite[Theorem 3.1]{RSV}, there exists $c_1=c_1(\alpha,d, M)
> 0$ such that
 \beq\label{e:gm}
\lim_{|x| \to 0} G^M (x) \, |x|^{d-2} \, =\, c_1.
 \eeq
Let $L_1=(2/c_1)^{2/(d-2)} \vee 2$ and take $ 0<\delta_1<c_1 (
1-L^{-(d-2)/2}). $
Using \eqref{e:gm}, we can choose a positive constant $R_1 >0$
 such that
\begin{equation}\label{e:3.4}
(c_1-\delta_1)\frac1{|x|^{d-2}}\, \le\, G^M(x) \qquad \hbox{when }  |x|\le
R_1.
\end{equation}
Thus, by Lemma \ref{t:Gorigin}, for every $|x| < R_1$
\begin{eqnarray*}
G^a(x) \ge G^M (x) \ge (c_1-\delta_1)\frac1{|x|^{d-2}}
\ge  L_1^{\frac{d -2}2}
\frac{c_1}{|L_1 x|^{d-2}}\,\ge \,2\, G^a(L_1x).
\end{eqnarray*}
\qed

The next proposition gives the interior estimates for $G_D^a$.

\begin{prop}\label{GI}
Suppose that $M>0$.
For any bounded and connected $C^{1,1}$ open set $D$ in $\bR^d$
there exists a positive constant $C_4$ such that for every $a
\in [0,M]$
 \beq\label{ub} G^a_D(x,y) \le C_4 \frac1{|x-y|^{d-2}} \qquad
\text{ for all }x,y \in D
 \eeq
 and
  \beq \label{lb} G^a_D(x,y) \,\ge\, C_4^{-1}\,
\frac1{|x-y|^{d-2}}\ \qquad \text{ when }2|x-y| \le  \delta_D(x)
\wedge \delta_D(y)  .
 \eeq
\end{prop}

\pf Since $G^a_D(x,y) \le G^a(x,y)$, \eqref{ub} is an immediate
consequence of Lemma \ref{t:Gorigin}. So we only need to show
\eqref{lb}. Without loss of generality we assume that $\delta_D(y) \le
\delta_D(x)$.

Recall that $L_1\geq  2$ and $R_1$ are the constants from Lemma
\ref{G_G}. By Lemmas \ref{t:Gorigin}--\ref{G_G} and \eqref{e:3.4},
we have
 \beq\label{e:lb}
G^a(x,y) - G^a(L_1x, L_1y) \,
\ge\, \frac12 G^a(x,y)\ge\, \frac12 G^M(x,y)
\,\ge\, c_1
\frac1{|x-y|^{d-2}},  \qquad \text{ when } |x-y| \le R_1
 \eeq
for some positive constant $c_1$.

\noindent {\it Case 1}: $L_1|x-y| \le   \delta_D(y)$.
We consider three subcases separately:

(a) $\delta_D(y) \le R_1$. Note that, since $L_1|x-y| \le \delta_D(y)$,
$$
|X^a_{\tau^a_{B(y, \delta_D(y))}}-y| \,\ge \, \delta_D(y) \, \ge
\,L_1 |x-y| .
$$
Thus by the fact that $G^a(\cdot)$ is radially decreasing and
\eqref{e:lb},
\begin{eqnarray*}
G^a_D(x,y) &\ge & G^a_{B(y, \delta_D(y))}(x,y)
\,=\, G^a(x,y) - \E_x\left[ G^a( X^a_{\tau^a_{B(y, \delta_D(y))}},y) \right]\\
&\ge & G^a(x,y) - G^a(L_1x, L_1y) \,\ge\, c_1 \frac1{|x-y|^{d-2}}.
\end{eqnarray*}

(b) $\delta_D(y) > R_1$ and $L_1|x-y| \le R_1$. In this case,
$|X^a_{\tau^a_{B(y, R_1)}}-y| \,\ge \, R_1 \,\ge \,L_1|x-y|$ and,
again
by the fact that
$G^a(\,\cdot\,)$ is radially decreasing and
\eqref{e:lb},
\begin{eqnarray*}
 G^a_D(x,y) &\ge & G^a_{B(y, R_1)}(x,y)
\,=\, G^a(x,y) - \E_x\left[     G^a( X^a_{\tau^a_{B(y, R_1)}}, y) \right]\\
 &\ge & G^a(x,y) - G^a(L_1x, L_1y)\,\ge\, c_1 \frac1{|x-y|^{d-2}}.
\end{eqnarray*}

(c) $\delta_D(y) > R_1$ and $L_1|x-y| > R_1$. In this case, we have
$ \delta_D(x) \ge \delta_D(y) \ge L_1|x-y|\ge  R_1$. Choose a point
$w \in
\partial B(y , R_1/(2L_1))$. Then from the argument in (b), we get
$$
G^a_D(w,y) \,\ge\, c_1 \frac1{(R_1/(2L_1))^{d-2}}.
$$
Since $D$ is a bounded and connected $C^{1,1}$ open set and $|x-w|
\le |x-y|+|y-w| \le \delta_D(y)/L_1 + R_1/(2L_1) \le
3\delta_D(y)/(2L_1)$, by Proposition \ref{uhp} and a chain argument,
we have
\begin{eqnarray*}
G^a_D(x,y) \,\ge\, c_2 G^a_D(w,y) \,\ge\, c_3
\frac1{(R_1/(2L_1))^{d-2}} &\geq & c_3 \, 2^{d-2}
(2L_1/R_1)^{2(d-2)} \, \frac1{|x-y|^{d-2}}.
\end{eqnarray*}

\noindent {\it Case 2}: $2|x-y| \le  \delta_D(y)< L_1 |x-y|$. Take  $x_0
\in \partial  B(y, \delta_D(y)/(L_1+1))$. Then
$$
|x-y| \le \frac{1}{2}  \delta_D(y) \leq L_1 |x_0-y|
= \frac{L_1}{L_1+1} \delta_D(y) \le   \delta_D(y)
\wedge \delta_D (x_0).
$$
Since
$D$ is a bounded and connected $C^{1,1}$ open set and
$|x_0-x| \le |x_0-y|+|y-x|\leq (\frac{1}{L_1+1}  + \frac12)
\delta_D(y)$, by Proposition \ref{uhp}, a chain argument and the
argument in the first case, there are constants $c_i=c_i(D, \alpha,
L_1, M)>0, i=4,5,6,$ such that
$$
G^a_D(x, y) \ge  c_4 \, G^a_D(x_0, y) \ge  c_5
\frac1{|x_0-y|^{d-2}}\, \ge  c_6 \frac{1}{|x-y|^{d-2}} .
$$
This completes the proof of the proposition.  \qed

Suppose that $D$ is a bounded and connected $C^{1,1}$ open set in
$\R^d$ with characteristics $(R,\Lambda)$ and corresponding
$\kappa$. Fix $z_0\in D$ with $\kappa R<\delta_D(z_0)<R$, and let
$\eps_1:=\kappa R/24$. For $x,y\in D$,  define $r(x,y): =
\delta_D(x) \vee\delta_D(y)\vee |x-y|$ and
$$
\sB(x,y):=\left\{z \in D:\, \delta_D(z) > \frac{\kappa}{2}r(x,y), \,
|x-z|\vee |y-z| < 5 r(x,y)  \right\}
$$
if $r(x,y) <\eps_1 $, and $\sB(x,y):=\{z_0 \}$ otherwise.

Put $C_5:=C_4 2^{d-2} \delta_D(z_0)^{-d+2}$. Then by \eqref{ub},
$$
G^a_D(\cdot, z_0)\leq C_5 \qquad \hbox{on } D \setminus B(z_0,
\delta_D(z_0)/2).
$$
Now we define
$$
g^a(x ):=  G^a_D(x, z_0) \wedge C_5.
$$
Note that if $\delta_D(z) \le 6 \eps_1$, then $|z-z_0| \ge
\delta_D(z_0) - 6 \eps_1 > \delta_D(z_0) /2$ since
$6\eps_1<\delta_D(z_0)/4$, and therefore $g^a(z )= G^a_D(z, z_0)$.

Using the uniform Harnack principle (Proposition \ref{uhp}) and
Proposition \ref{l:Green_L}, the following form of Green function
estimates follows from \cite[Theorem 2.4]{H}.
\medskip

\begin{thm}\label{t:Gest}
Suppose that $M>0$. For any bounded and connected $C^{1, 1}$ open
set $D$ in $\bR^d$, there exists $C_6=C_6(D, M, \alpha)>0$ such that
for all $x, y$ in $D$ and all $a \in (0,M]$
$$
C_6^{-1}\,\frac{g^a(x) g^a(y)}{g^a(A)^2} \,|x-y|^{-d+2} \,\le\,
G^a_D(x,y) \,\le\, C_6\,\frac{g^a(x) g^a(y)}{g^a(A)^2}\,
|x-y|^{-d+2},
$$
where $A\in \sB(x, y)$.
\end{thm}

Suppose  $D$ is a bounded and connected $C^{1, 1}$ open set. For
all $x,y \in D$, we let $Q_x$ and $Q_y$ be points on $\partial D$
such that $\delta_D(x)=|x-Q_x|$ and $\delta_D(y)=|y-Q_y|$
respectively. It is easy to check that  if $r(x,y) < \eps_1$
\begin{equation}\label{e:AinB}
A_{r(x,y)}(Q_x),\, A_{r(x,y)}(Q_y) \,\in\, \sB(x,y).
\end{equation}
(Recall that, for any $Q\in \partial D$, $A_r(Q)$ is a point such
that $B(A_r(Q),\kappa r)\subset D\cap B(Q, r)$.)
Indeed, by the definition of $ A_{r(x,y)}(Q_x)$,
$\delta_D(A_{r(x,y)}(Q_x)) \ge \kappa r(x,y) > \kappa r(x,y) /2$.
Moreover,
$$
|x- A_{r(x,y)}(Q_x)|\, \le\, | x-Q_x|+|Q_x -  A_{r(x,y)}(Q_x)|
\,\le\, \delta_D(x) + r(x,y) \,\le\, 2 r(x,y)
$$
and $ |y- A_{r(x,y)}(Q_x)| \le |y-x| +|x- A_{r(x,y)}(Q_x)| \le 3
r(x,y) $. This verifies the claim \eqref{e:AinB}.

Recall the fact that  $g^a(z )=  G^a_D(z, z_0)$ if $\delta_D(z) < 6
\eps_1$. By Theorem \ref{t:main} and the fact that $c_0^{-1} r(x,y)
\le \delta_D(A_{r(x,y)}(Q_y)) \le c_0 r(x,y)$ for some $c_0>1$,
there exists $c_1>1$ such that for every $a \in (0,M]$ and all $x ,y
\in D$ with $\delta_D(x) < 6 \eps_1$ and $\delta_D(y) < 6 \eps_1$,
\begin{equation}\label{e:com}
c_1^{-1}  \frac{\delta_D(x)}{r(x,y)} \le  \frac{g^a(x)} {
g^a(A_{r(x,y)}(Q_x))} = \frac{G^a_D(x, z_0)} {
G^a_D(A_{r(x,y)}(Q_x), z_0)} \le c_1 \frac{\delta_D(x)}{r(x,y)}
\end{equation}
and
\begin{equation}\label{e:com2}
c_1^{-1}  \frac{\delta_D(y)}{r(x,y)} \le  \frac{g^a(y)} {
g^a(A_{r(x,y)}(Q_y))} = \frac{G^a_D(y, z_0)} {
G^a_D(A_{r(x,y)}(Q_y), z_0)} \le c_1 \frac{\delta_D(y)}{r(x,y)}.
\end{equation}

\bigskip

\noindent {\bf Proof of Theorem \ref{t-main-green} when $d\geq 3$:}
First we assume that $D$ is connected.

Combining inequalities \eqref{e:com} and \eqref{e:com2} with
Proposition \ref{GI}, Theorem \ref{t:Gest} and the fact that
\begin{equation}\label{e:gg2}
\frac{ \delta_D(x) \delta_D(y)}{(r(x,y))^2} \le  \left(1\wedge
\frac{  \delta_D(x) \delta_D(y)}{ |x-y|^{2}}\right) \le  \frac94
\frac{ \delta_D(x) \delta_D(y)}{(r(x,y))^2}
\end{equation}
(see \cite{Bo}), we get the inequalities \eqref{e:1.3}.

Next we assume that $D$ is not connected. Let $(R, \Lambda)$ be the
$C^{1,1}$ characteristics of $D$. Note that $D$ has only finitely
many components and the distance between any two distinct components
of $D$ is at least $R>0$. Assume first that $x$ and $ y$ are  in
two  distinct  components of $D$. Let $D(x)$ be the component of $D$
that contains $x$. Then by the strong Markov property and the L\'evy
system \eqref{e:levy} of $X^a$, we have
$$
G_D^a(x, y)= \E_x \left[ G^a_D\big(X^a_{\tau^a_{  D(x)}}, y\big) \right] =
\E_x \left[ \int_0^{\tau^a_{ D(x)}} \left( \int_{D\setminus D(x)}
j^a (|X_s^a-z|) G^a_D(z, y) dz \right)  ds \right].
$$
Consequently,
\begin{equation}\label{e:3.12}
j^a({\rm diam} (D)) \, \E_x[\tau^a_{D(x)}] \, \int_{D\setminus D(x)}
G^a_D(y, z) dz  \leq G^a_D(x, y) \leq j^a(R) \, \E_x[\tau^a_{
D(x)}] \, \int_{D\setminus D(x)} G^a_D(y, z) dz.
\end{equation}
Applying the two-sided estimates \eqref{e:1.3} established in the
first part of this proof to $D(x)$, we get
\begin{equation}\label{e:3.13}
c_1^{-1} \delta_D(x)=c_1^{-1} \delta_{D(x)}(x) \le \E_x
\left[ \tau^a_{D(x)} \right] \le c_1\delta_{D(x)}(x) =c_1 \delta_D(x)
\end{equation}
for some $c_1=c_1(D, M,\alpha)>1$.
Clearly,
using \eqref{e:3.13},
$$
\int_{D\setminus D(x)} G^a_D(y, z) dz \geq \int_{D(y)} G^a_{D(y)}
(y, z) dz =\E_y [\tau^a_{D(y)}] \geq c_2 \, \delta_D(y).
$$
On the other hand, it follows from \eqref{ub} that
$ \sup_{z\in D, a\in (0, M]}\E_z [\tau^a_D] \leq c_3<\infty$.
Moreover by \eqref{e:3.13} and
the L\'evy system \eqref{e:levy} of $X^a$,
\begin{eqnarray*}
\int_{D\setminus D(x)} G^a_D(y, z) dz
 &\leq & \E_y \big[
\tau^a_D\big]  = \E_y \Big[ \tau^a_{D(y)}\Big]
+ \E_y \Big[\E_{X_{\tau^a_{D(y)}}} [\tau^a_D ]\Big]\\
&\leq & c_4 \, \delta_D(y) +
\E_y \left[ \int_0^{\tau^a_{D(y)}} \left(
\int_{D\setminus D(y)}
j^a (|X_s- z|) \, \E_z [\tau^a_D] dz \right)  ds \right] \\
&\leq & c_4 \, \delta_D(y) + c_5
 j^M(R) \E_y \big[ \tau^a_{D(y)}\big]
\leq c_6 \, \delta_D(y).
\end{eqnarray*}
We conclude from the last three displays, \eqref{e:3.12} and the
form of $j^a$ given in \eqref{e:j} that there is a constant
$c_7=c_7(D, M, \alpha)\geq 1$ such that for every $a\in (0, M]$,
\begin{equation}\label{e:dfc}
c_7^{-1} a^\alpha \, \delta_D(x) \delta_D(y) \leq G_D(x, y) \leq
 c_7  a^\alpha \, \delta_D(x) \delta_D(y).
\end{equation}
Since for $x$ and $ y$ in different components of $D$, $R\leq
 |x-y| \leq {\rm diam}(D)$, we have established  \eqref{e:1.3}.

Now we assume that $x, y$ are in the same component $U$ of $D$.
Applying \eqref{e:1.3} to $U$ we get
$$G^a_D(x,y)\ge G^a_U(x,y) \ge  c_8 \frac{1} {|x-y|^{d-2}} \left(1\wedge
\frac{  \delta_U(x) \delta_U(y)}{ |x-y|^{2}}\right)=
c_8 \frac{1} {|x-y|^{d-2}} \left(1\wedge\frac{  \delta_D(x) \delta_D(y)}{ |x-y|^{2}}\right).
$$
For the upper bound, we use the strong Markov property, the L\'evy
system \eqref{e:levy}, and \eqref{e:3.13}--\eqref{e:dfc} to get
\begin{eqnarray}
&&G^a_D(x,y)\,=\, G^a_U(x,y) +\E_x \left[ G^a_D(
X^a_{\tau^a_{  U}}, y) \right] \nonumber\\
&\le& c_9  \frac{1} {|x-y|^{d-2}} \left(1\wedge
\frac{  \delta_D(x) \delta_D(y)}{ |x-y|^{2}}\right) +
\E_x \left[ \int_0^{\tau^a_{ U}} \left( \int_{D\setminus U}
j^a (|X_s^a-z|) G^a_D(z, y) dz \right)  ds \right] \nonumber\\
&\le& c_9  \frac{1} {|x-y|^{d-2}} \left(1\wedge
\frac{  \delta_D(x) \delta_D(y)}{ |x-y|^{2}}\right) +
j^M(R) \, \E_x[\tau^a_{
U}] \, \int_{D\setminus U} G^a_D(y, z) dz \nonumber\\
&\le& c_9  \frac{1} {|x-y|^{d-2}} \left(1\wedge \frac{  \delta_D(x)
\delta_D(y)}{ |x-y|^{2}}\right) +  c_{10} \delta_D(x) \delta_D(y) \,
\int_{D\setminus U} \delta_D(z) dz. \label{e:fip}
\end{eqnarray}
Since the boundedness of $D$ implies
$$
\delta_D(x) \delta_D(y)  \le
c_{11} \frac{1} {|x-y|^{d-2}} \left(1\wedge
\frac{  \delta_D(x) \delta_D(y)}{ |x-y|^{2}}\right),
$$
we have from \eqref{e:fip}
$$
G^a_D(x,y) \le
c_{12}\frac{1} {|x-y|^{d-2}} \left(1\wedge
\frac{  \delta_D(x) \delta_D(y)}{ |x-y|^{2}}\right).
$$
\qed

Define
\begin{equation}\label{e:a}
a (x, y,z,w) := \begin{cases} 1
\quad &\hbox{when } x, y, z, w \hbox{ are in the same component of } D, \\
 a^{-\alpha} \quad &\hbox{when } x, y \in D(x),
  z\notin D(x) \hbox{ and } w\in D(z), \\
a^\alpha
\quad &\hbox{when }  x, w \in D(x), \#(\{y,z\} \cap D(x))=1, \\
a^\alpha\quad &\hbox {when } x,y,z,w \hbox{ are all in different components of } D,\\
 a^{2\alpha}\quad &\hbox{when }   x, w \in D(x),\{y,z\}
 \cap D(x) =\emptyset.
\end{cases}
\end{equation}

The next theorem will be used in Section \ref{s:pr}.

\begin{thm}[Generalized 3G theorem]\label{t:3Gt}
Suppose that $M>0$. For any bounded $C^{1, 1}$ open set $D$ in
$\bR^d$, there exists a constant $C_8=C_8(D, \alpha, M)$ such that
for all $x, y, z, w \in D$ and $a \in (0, M]$,
\begin{align}
&\frac{G^a_D(x,y) G^a_{D}(z,w)} { G^a_{D}(x,w)} \nonumber
\\&
 \le C_8 a(x,y,z,w)
\left(\frac{|x-w|\wedge |y-z|}{|x-y|  } \vee 1 \right)
\left(\frac{|x-w|\wedge |y-z|}{|z-w|} \vee 1 \right)
\frac{|x-w|^{d-2}} {|x-y|^{d-2} |z-w|^{d-2}}.\label{3G}
\end{align}
\end{thm}

 \pf
Recall that  $r(x,y)=
\delta_D(x) \vee\delta_D(y)\vee |x-y|$ and
 let
 $$
g_D(x, y) :=  \frac{1} {|x-y|^{d-2}} \frac{ \delta_D(x) \delta_D(y)}{(r(x,y))^2} \quad
\mbox{and} \quad H(x,y,z,w)\,:=\,\frac{ |x-w|^{d-2}}
{|x-y|^{d-2} |z-w|^{d-2}   }.
$$
 By Theorem \ref{t-main-green} for the case $d\geq 3$ and \eqref{e:gg2},
\begin{align}
\frac{G^a_D(x,y) G^a_{D}(z,w)} { G^a_{D}(x,w)} &\le c_1 a(x,y,z,w) \frac{g_D(x,y) g_{D}(z,w)} { g_{D}(x,w)}
\label{e:gtt0}\\ &= c_1 a(x,y,z,w)
\frac{\delta_D(y)\delta_D(z) r(x,w)^2}{r(x,y)^2 r(z,w)^2}
H(x,y,z,w).\label{e:gtt}
\end{align}

\begin{enumerate}
\item If
$|x-w| \le  \delta_D(x) \wedge \delta_D(w)$, $ g_D(x,w)  \ge
|x-w|^{-d+2}. $ Thus, by \eqref{e:gtt0}
$$
\frac{G^a_D(x,y) G^a_{D}(z,w)} { G^a_{D}(x,w)} \le  c_1 a(x,y,z,w)
H(x,y,z,w).
$$

 \item
Note that if $y=z$, since $r(x,w) \le 2 r(x,y) + 2 r(y,w)$,
 \begin{align*}
\frac{\delta_D(y)\delta_D(y) r(x,w)^2}{r(x,y)^2 r(y,w)^2} \le 8
\left(         \frac{\delta_D(y)^2}{r(y,w)^2} +\frac{\delta_D(y)^2 }{r(x,y)^2 }  \right)\le 8.
\end{align*}
 Thus
 \begin{align}
\frac{g_D(x,y) g_{D}(y,w)} {g_{D}(x,w)} \le 8
H(x,y,y,w).\label{e:gtt1}
 \end{align}

Now consider the case $|y-z| \leq \delta_D(y) \wedge \delta_D(z)$.  In
this case $ g_D(y,z)  \ge  |y-z|^{-d+\alpha}. $ Thus, using
\eqref{e:gtt1}, we obtain that
\begin{eqnarray}
&& \frac{g_D(x,y) g_{D}(z,w)} { g_{D}(x,w)}\,=\,
\frac{g_D(x,y) g_{D}(y,z)} { g_{D}(x,z)}  \frac{g_D(x,z) g_{D}(z,w)}
{ g_{D}(x,w)}\frac{1}{ g_{D}(y,z) }\nonumber\\
&&\le \,64 \,\frac{|x-z|^{d-2}} {|x-y|^{d-2}
|y-z|^{d-2}}\frac{|x-w|^{d-2}}
{|x-z|^{d-2} |z-w|^{d-2}}\frac{1}{   g_{D}(y,z)   }\nonumber\\
&&=\, 64 \, \frac{ |x-w|^{d-2}} {|x-y|^{d-2}
|y-z|^{d-2}|z-w|^{d-2} }\frac{1}{ g_{D}(y,z)}.
\label{3G_est53}
\end{eqnarray}
Thus, by \eqref{e:gtt0} and \eqref{3G_est53}, we have
$$ \frac{G^a_D(x,y) G^a_{D}(z,w)} { G^a_{D}(x,w)} \le c_2 a(x,y,z,w) H(x,y,z,w).
$$
 \item
Now we assume that $|x-w| >  \delta_D(x) \wedge \delta_D(w)$ and $|y-z|
> \delta_D(y) \wedge\delta_D(z)$. Since $\delta_D(x) \vee \delta_D(w) \le \delta_D(x) \wedge
\delta_D(w) +|x-w|$, using the assumption $ \delta_D(x) \wedge
\delta_D(w)< |x-w|$, we obtain $r(x,w) <2|x-w|$. Similarly,
$r(y,z)<2|y-z|$. By \eqref{e:gtt}, we only need to show that
\begin{equation}\label{e:33}
\frac{\delta_D(y)\delta_D(z) r(x,w)^2}{r(x,y)^2 r(z,w)^2} \le
c_3\left(\frac{|x-w|\wedge |y-z|}{|x-y|  } \vee 1 \right)
\left(\frac{|x-w|\wedge |y-z|}{|z-w|} \vee 1 \right).
\end{equation}
Since
$r(x,w) \le 2 r(x,y) + 2 r(y,w) \le 2r(x,y) + 4r(y, z)+ 4 r(z,w)$,
we
have
\begin{eqnarray*}
\frac{\delta_D(y)\delta_D(z) r(x,w)^2}{r(x,y)^2 r(z,w)^2} &\le& c_4
\left(         \frac{\delta_D(y)\delta_D(z)}{r(z,w)^2} +\frac{\delta_D(y)
\delta_D(z) }{r(x,y)^2 } +\frac{\delta_D(y)\delta_D(z) r(y,z)^2}
{r(x,y)^2 r(z,w)^2}             \right)\\
 &\le& c_4  \left(         \frac{\delta_D(y)}{r(z,w)} +
 \frac{\delta_D(z) }{r(x,y) } +\frac{ r(y,z)^2}{r(x,y) r(z,w)}
 \right)\\
 &\le& c_4  \left(         \frac{r(y,z)}{r(z,w)} +
 \frac{r(y,z) }{r(x,y) } +\frac{ r(y,z)^2}{r(x,y) r(z,w)}
 \right),
\end{eqnarray*}
which is, by \cite[Lemma 3.15]{KL}, less than or equal to
$$
2c_4\left(\frac{r(y,z) }{r(x,y) } \vee 1 \right)\left(\frac{r(y,z)}{r(z,w)}  \vee 1 \right).
$$
On the other hand, clearly
\begin{eqnarray*}
\frac{\delta_D(y)\delta_D(z) r(x,w)^2}{r(x,y)^2 r(z,w)^2} &=&
\frac{\delta_D(y)\delta_D(z)}{r(x,y) r(z,w)} \frac{r(x,w)^2}{r(x,y)^2 r(z,w)^2} \\
&\le & \left(\frac{r(x,w) }{r(x,y) } \vee 1 \right)
\left(\frac{r(x,w)}{r(z,w)}  \vee 1 \right).
\end{eqnarray*}
Thus
$$
\frac{\delta_D(y)\delta_D(z) r(x,w)^2}{r(x,y)^2 r(z,w)^2}
\le\, c_{5} \left(\frac{r(y,z) \wedge r(x,w)}{r(x,y)} \vee 1
\right) \left(\frac{r(y,z)\wedge r(x,w)}{r(z,w)} \vee 1
\right) .
$$
Now applying the fact that $r(x,w)<2|x-w|$,
$r(y,z)<2|y-z|$, $r(x,y) \ge |x-y|$ and $r(z,w) \ge
|z-w|$,  we arrive at (\ref{e:33}).
\end{enumerate}
We have proved the theorem.
 \qed

Note that, since we consider disconnected open sets too, we can not
apply \cite[Theorem 1.1]{KL} directly to get the generalized 3G
theorem.

\medskip

Taking $y=z$ in (\ref{3G}), we get  the classical 3G estimates, that
is,
$$
\frac{G^a_D(x,z) G^a_D(z,w)} { G^a_D(x,w)} \, \le\, \,C_8\,
a(x, z, z, w)\frac{|x-w|^{d-2}} {|x-z|^{d-2} |z-w|^{d-2}} \,\le\,
C_8(M^{2\alpha} \vee 1)\frac{|x-w|^{d-2}} {|x-z|^{d-2} |z-w|^{d-2}}.
$$

\section{Two Dimensional Case}

In this section we assume $d=2$ and prove Theorem \ref{t-main-green}
for this case. Unlike the case of $d\geq 3$, due to the recurrence
of planar Brownian motions, the methods in \cite{Bo1, H} are not
applicable in dimension $d=2$ even though we have the Harnack and
boundary Harnack principles. We use a capacitary approach and some
recent results on the subordinate killed Brownian motions instead.

First we derive the lower bound. The method we use relies on
comparing the process $X^{a,D}$, which is the killed subordinate
Brownian motion, with another process, the subordinate killed
Brownian motion. This method also works for dimensions $d\ge 3$.

To be more precise, let $D$ be a bounded open set in $\bR^2$ and
$X^{0, D}$ the killed Brownian motion in $D$. Let $(T^a_t: t\ge 0)$
be a subordinator independent of $X^0$ which can be written as
$T^a_t=t+a^2T_t$ where $(T_t: t\ge 0)$ is an $\alpha/2$-stable
subordinator. The process $(Z^{a, D}_t: t\ge 0)$ defined by $Z^{a,
D}_t=X^{0, D}_{T^a_t}$ is called a subordinate killed Brownian
motion in $D$. Let $u^a$ be the potential density of $T^a$ (see
\eqref{e:pd4sub}). It follows from \cite{SV06} that the Green
function $R^a_D(x, y)$ of $Z^{a, D}$ is given by
\begin{equation}\label{gfn4sks}
R^a_D(x, y)=\int^{\infty}_0 p^0_D(t, x, y)u^a(t)dt,
\end{equation}
where $p^0_D(t, x, y)$ is the transition density of the killed
Brownian motion $X^{0, D}$.
It is well known
(see, for instance, \cite[Proposition 3.1]{SV08})
that
\begin{equation}\label{gfcrln}
R^a_D(x, y)\le G^a_D(x, y), \qquad (x, y)\in D\times D.
\end{equation}

\begin{thm}\label{T:gfcesks}
Suppose that $M>0$. For any bounded $C^{1, 1}$ \emph{connected} open
set $D$ in $\R^2$, there exists a positive constant $C_9=C_9(\alpha,
M, D)$ such that for all $x,y \in D$ and all $a\in (0, M]$,
$$
G^a_D(x,y) \ge R^a_D(x, y)\ge C_9\,
\log
\left(1+\frac{\delta_D(x)\delta_D(y)}{|x-y|^2} \right).
$$
\end{thm}

\pf First recall the following lower bound for the transition
density of the killed Brownian motion $X^{0, D}$ obtained in
\cite{Zq3} which states that for any $A>0$, there exist positive
constants $c_0$ and $c_1$ such that for any $t\in (0,A]$ and any
$x,y\in D$,
\begin{equation}\label{lower bound for p}
p_D^0(t,x,y)\ge c_0\left(1\wedge \frac{\delta_D(x)\delta_D(y)}{t}\right) \,
t^{-1}\exp\left(-\frac{c_1 |x-y|^2}{t}\right)\, .
\end{equation}
It follows from \eqref{e:pd4sub} that
\begin{equation}\label{e:relpd}
u^a(t)=u^1(a^{\frac{2\alpha}{2-\alpha}}t) \qquad \hbox{for } t>0.
\end{equation}
Let $T=({\rm diam}(D))^2$. Since $u^1(t)$ is a completely monotone
function with $u^1(0+)=1$, by \eqref{e:relpd},
 for any $a\in (0,M]$
\begin{align}
u^a(t)\ge  u^1(M^{\frac{2\alpha}{2-\alpha}}T), \qquad t\in (0, T].
\label{e:stlb4pd}
\end{align}
By a change of variables $s= \frac{|x-y|^2}{t}$, we have
\begin{eqnarray}
\int_0^T\left(1 \wedge \frac{\delta_D(x)\delta_D(y)}{t}\right)
t^{-1}e^{-c_1\frac{|x-y|^{2}}{t}} dt
=\int_{\frac{|x-y|^2}{T}}^\infty \left(1\wedge \frac{   \delta_D(x)
\delta_D(y)\, s}{ |x-y|^2 }\right) s^{-1}   e^{-c_1s}ds.
\label{ID:5.1}
\end{eqnarray}
Define
\begin{equation}\label{e:f1}
f_D(x, y)= \frac{\delta_D(x) \delta_D(y)}{|x-y|^2}.
\end{equation}
Since
$1/f_D(x,y) \geq |x-y|^2/{\rm diam}(D)^{2}=|x-y|^2/T$,
we split the last integral into two parts:
\begin{eqnarray}
\int_{1}^\infty \left(1\wedge \frac{   \delta_D(x) \delta_D(y)\, s}{
|x-y|^2 }\right)s^{-1}e^{-c_1s} ds \geq  \left(  \int_{1}^\infty
s^{-1}\,  e^{-c_1s} ds \right)   \left(1\wedge \frac{  \delta_D(x)
\delta_D(y)}{ |x-y|^2 }\right)  , \label{e:ch1}
\end{eqnarray}
and
\begin{eqnarray}
\int_{\frac{|x-y|^2}{T}}^{1 }s^{-1} \left(1\wedge
(f_D(x,y)\, s)\right) ds
&\ge &\int_{\frac{|x-y|^2}{T}}^{1 } s^{-1} {\bf 1}_{\{s\geq
1/f_D(x,y)\}} ds \nonumber \\
&=& \log (f_D(x,y) \vee 1).  \label{e:ch7}
\end{eqnarray}
Combining \eqref{lower bound for p} and
\eqref{ID:5.1}--\eqref{e:ch7}, we have
\begin{eqnarray*}
\int_0^T p^0_D(t, x, y)dt&\ge &c_2(1\wedge f_D(x,y)) + c_2\log (f_D(x,y) \vee 1)\\
&\ge &  c_3 \log \left( 1+ \frac{  \delta_D(x) \delta_D
(y)}{|x-y|^{2} }\right) .
\end{eqnarray*}
So it follows from \eqref{gfn4sks}--\eqref{gfcrln} and
\eqref{e:stlb4pd} that
$$
G^a_D(x,y)\ge R^a_D(x, y)\,\ge \,u^1(M^{\frac{2}{2-\alpha}}T)
\int_0^Tp^0_D(t, x, y) dt \, \ge\, c_4\log
\left(1+\frac{\delta_D(x)\delta_D(y)}{|x-y|^2} \right).
$$
\qed

Integrating the estimate in Theorem \ref{T:gfcesks}
with respect to $y$ yields the following corollary.

\begin{cor}\label{c-lb-exit}
Suppose that $M>0$. For any bounded
connected $C^{1, 1}$
open
set $D$ in $\R^2$, there exists a positive constant
$C_{10}=C_{10}(\alpha, M, D)$ such that for all $x\in D$ and all
$a\in (0, M]$,
$$
\E_x [\tau_D^a]\ge C_{10}\delta_D(x)\, .
$$
\end{cor}

The inequalities in the next lemma can be proved by elementary
calculus and will be used
several times without being mentioned
explicitly.

\begin{lemma}\label{l:Gest0}
For any $L>0$, there exists a constant $C_{11}=C_{11}(L)>1$ such
that
$$
C_{11}^{-1} b \le \log(1+b) \le b  \quad \text{for any }0<b\le L
$$
and
$$
C_{11}^{-1} \log(1+s)  \le \log(1+Ls) \le C_{11} \log(1+s)  \quad
\text{ for any }0<s < \infty.
$$
\end{lemma}

Using Corollary \ref{c-lb-exit}, Theorem \ref{T:gfcesks} can be extended to
general (not necessarily connected) bounded $C^{1,1}$ open sets.
Recall that $g_D^a$ is defined by \eqref{e:1.2}.

\begin{thm}\label{T:gfcesks2}
Suppose that $D$ is a bounded $C^{1, 1}$ open set in $\bR^2$ with
characteristics $(R, \Lambda)$. There exists a positive constant
$C_{12}=C_{12}(\alpha, M, D)$ such that for all $x,y \in D$ and all
$a\in (0, M]$,
$$
G^a_D(x,y) \ge  C_{12}\,
 g_D^a (x, y).
$$
\end{thm}

\pf
Recall that $f_D(x, y)$ is defined in \eqref{e:f1}.
If $x$ and $y$ are in the same component,
 say $x,y\in U$, then by monotonicity,
\begin{equation}\label{lower bound same}
G^a_D(x,y)\ge G^a_{\, U}(x,y)
\ge c_1 \log(1+f_{U}(x, y)) = c_1 \log(1+f_D(x, y)).
\end{equation}
If $x,y$ are in the different components of $D$, using Corollary
\ref{c-lb-exit} and Lemma \ref{l:Gest0}, and by following the second
part of the proof of Theorem \ref{t-main-green} in case $d\geq 3$
(that is, the paragraph containing
\eqref{e:3.12}-\eqref{e:3.13}), we get
\begin{eqnarray*}
G^a_D(x,y)\ge   c_2 a^\alpha
\delta_D(x)\delta_D(y)  \ge   c_2R^2 a^\alpha
\frac{\delta_D(x)\delta_D(y)}{|x-y|^2}
 \ge
 c_{3} a^\alpha
\log \left(1+ \frac{\delta_D(x)\delta_D(y)}{|x-y|^2}\right)\, .
\end{eqnarray*}
This completes the proof of the theorem.
\qed

Recall that when $d\ge 2$ and $a>0$, any non-empty open set
$D\subset\bR^d$ is Greenian for $X^a$. For any Greenian open set
$D$, any  Borel subset $A$  of $D$ and $a\geq 0$, we define
\begin{eqnarray}
\mbox{Cap}^a_{D}(A)&:=&\sup \Big\{\eta(A): \eta
\mbox{ is a measure supported on }A \nonumber \\
&&\hskip 1truein \mbox{ with } \int_{D} G^{a}_D(x,y) \eta(dy) \le 1
\Big\}\, . \label{cap1}
\end{eqnarray}
It is known  (cf. \cite{FOT}) that for any open subset $A$ of $D$,
$$
\mbox{Cap}^a_{D}(A)\,=\,\inf \left\{ \sE^a(u, u): \ u \in
W^{1,2}(\R^d), u=0 \text{ on } D^c,  u\ge 1 \text{ a.e. on } A
\right\}
$$
and for any Borel subset $A$ of $D$,
$$
\mbox{Cap}^a_{D}(A)\,=\,\inf\{\mbox{Cap}^a_{D}(B) : A \subset B \text{
and } B \text{ is open}\}.
$$
Since $\sE^0\le \sE^a$, for any Greenian open set  $D \subset \bR^d$
and every $a \in [0,M]$
\begin{equation}\label{e:cap}
\mbox{Cap}^0_{D}(A)\le \mbox{Cap}^a_{D}(A)  \quad \text{for every }
A \subset D.
\end{equation}

\begin{lemma}\label{l:Cap1}
There exists $C_{13}>0$ such that
$$
\mathrm{Cap}^0_{B(0,1)}(\overline{B(0,r)}) \ge
\frac{C_{13}}{\log(1/r)}  \quad \text{for every } r \in (0, 3/4).
$$
\end{lemma}
\pf
Recall that (see, e.g., \cite[p. 178]{CZ})
\begin{equation}\label{e-green-ball}
G^0_{B(0,1)}(x,y)=\frac{1}{2\pi}\log\left(1+\frac{
(1-|x|^2)(1-|y|^2)}{|x-y|^2}\right)\, .
\end{equation}
Let $\mathcal{P}$ denote the family of all probability measures on
$\overline{B(0,r)}$. It follows from \cite[p.159]{Fu} that
\begin{equation}\label{e-formula-for-capacity}
\mathrm{Cap}^0_{
B(0,1)
}(
\overline{B(0,r)})=\left(\inf\limits_{\mu\in
\mathcal{P}}\sup\limits_{x\in \mathrm{supp}(\mu)}
G^0_{B(0,1)}\mu(x)\right)^{-1}\, .
\end{equation}
Let $m_r$ be the normalized Lebesgue measure on $\overline{B(0,r)}$.
By \eqref{e-formula-for-capacity},
\begin{equation}\label{e-estimate-for-capacity}
\mathrm{Cap}^0_{B(0,1)}(
\overline{B(0,r)})\ge \frac{1}{\sup\limits_{x \in
\overline{B(0, r)}} G^0_{B(0,1)}m_r(x)}\, .
\end{equation}
Further, by using symmetry in the first equality, and
\eqref{e-green-ball} in the second line, we have
\begin{eqnarray*}
&&\sup\limits_{x\in \overline{B(0,r)}} G^0_{B(0,1)}m_r(x)
=G^0_{B(0,1)}m_r(0)=\int_{B(0,r)}G^0_{B(0,1)}(0,y)\, m_r(dy)\\
&&=\frac{1}{\pi r^2}
\int_{B(0,r)}\frac{1}{2\pi}\log\frac{1}{|y|^2}\, dy\,=\,
\frac{1}{\pi r^2}\, \frac{r^2}{2}\left(1+2\log\frac1r\right)\,
\le\,
c \log\frac1r\, ,
\end{eqnarray*}
for some constant $c>0$. This together with
\eqref{e-estimate-for-capacity} yields the desired capacity
estimate. \qed

For any Borel subset $V$, we use $\sigma^a_V$ to denote the first
hitting time of $V$ by $X^a$: $\sigma^a_V=\inf\{t>0: X^a_t\in V\}$.

\begin{lemma}\label{l:Cap2}
Suppose that $M>0$. There exists $C_{14}>0$ such that for every $a
\in (0,M]$, any
Greenian
open set $D$ in $\bR^2$ containing $B(0,1)$
and any $x \in \overline{B(0, \frac34)}$
$$
G^a_{D}(x,0)\,\le\, \frac{C_{14}}{\mbox{\rm Cap}^0_{D}
\big(\overline{B(0, |x|/2)} \big)} \, \P_x\left(  \sigma^a_{\overline{B(0,
|x|/2 )}} <   \tau_{D}\right) .
$$
\end{lemma}

\pf Fix $x \in \overline{B(0, 3/4)}$ and let $r := |x|/2$. Since
$\overline{B(0, r )}$ is a compact subset of $D$, there exists a
capacitary measure $\mu^a_{r}$ for $\overline{B(0, r )}$ with
respect to $X^{a, D}$ such that
$$
\mbox{\rm Cap}^a_{D} (\overline{B(0, r )}) = \mu^a_{r}
(\overline{B(0, r )})
$$
(see, for example, \cite[Section VI.4]{BG} for details). Then by
Proposition \ref{uhp}, we have
\begin{eqnarray}
\int_{\overline{B(0, r )}}   G^a_{D}(x,y)  \mu^a_{r} (dy) &\ge&
\left(\inf_{y \in \overline{B(0, r )}}  G^a_{D}(x,y) \right)
\mu^a_{r}(\overline{B(0, r )})\nonumber \\
&\ge &c_1  G^a_{D}(x,0) \mbox{\rm Cap}^a_{D} (\overline{B(0, r )})\nonumber\\
&\ge &c_1  G^a_{D}(x,0) \mbox{\rm Cap}^0_{D} (\overline{B(0, r )})
\label{e:ddd}
\end{eqnarray}
for some constant $c_1>0$. In the last inequality above, we have
used (\ref{e:cap}).

On the other hand,
\begin{eqnarray}
\int_{\overline{B(0, r )}}   G^a_{D}(x,y)  \mu^a_{r} (dy) &=&
\int_{\overline{B(0, r )}} \E_x\left[  G^a_{D}(X^{a, D}_{
\sigma^a_{\overline{B(0, r )}}},y)\right]\
\mu^a_{r} (dy)\nonumber \\
&\le& \left(\sup_{w \in \overline{B(0, r )}}  \int_{\overline{B(0, r
)}}     G^a_{D}(w,y) \mu^a_{r} (dy) \right)
\P_x\left(  \sigma^a_{\overline{B(0, r )}} <   \tau_{D}\right)\nonumber\\
&\le& \P_x\left(  \sigma^a_{\overline{B(0, r )}} <   \tau_{D}\right)
\label{e:ddddd}.
\end{eqnarray}
In the last inequality above, we have used
\eqref{cap1}.

Combining (\ref{e:ddd})-(\ref{e:ddddd}) we have
$$
G^a_{D}(x,0) \le \frac{c_1^{-1}}{\mbox{\rm Cap}^0_{D}
\big(\overline{B(0, r )}\big)} \P_x\left(  \sigma^a_{\overline{B(0, r )}} <
\tau_{D}\right).
$$
\qed

\begin{cor}\label{c:Cap2}
Suppose that $M>0$. There exists $C_{15}>0$ such that for every $a
\in (0,M]$ and every $x \in \overline{B(0, 3/4)}$
$$
G^a_{B(0,1)}(x,0)\,\le\, C_{15}\,\log \left(1/|x| \right).
$$
\end{cor}

\pf It follows  from Lemmas \ref{l:Cap1}-\ref{l:Cap2} that
\begin{eqnarray*}
G^a_{B(0, 1)}(x,0) &\le& \frac{C_{14}}{\mbox{\rm Cap}^0_{B(0, 1)}
\big(\overline{B(0, |x|/2 )}\big)} \P_x\left(
\sigma^a_{\overline{B(0, |x|/2 )}} <
\tau_{B(0, 1)}\right)\nonumber\\
&\le& C_{14}C_{13}^{-1} \log \left(2/|x| \right) \le c \log
\left(1/|x| \right)
\end{eqnarray*}
for some constant $c>0$.
\qed

\begin{lemma}\label{l:Cap3}
Suppose that $M>0$ and that $D$ is a bounded $C^{1, 1}$ open set in
$\R^2$ with characteristics $(R, \Lambda)$. There exists
$C_{16}=C_{16}(D)>0$ such that for every $a \in (0,M]$ and all $x,y
\in D$ with $|x-y| \le \frac34 \delta_D(x) < \frac34 R$,
$$
G^a_{D}(x,y)\,\le\, C_{16}\,\log \left(\frac{\delta_D(x)}{|x-y|}
\right).
$$
\end{lemma}

\pf By our assumption, $D$ satisfies the uniform exterior ball
condition with radius $R>0$.

Fix $x,y \in D$ with $|x-y| \le \frac34 \delta_D(x)$ and let
$r:=\delta_D(x)$. Since $r < R$, without loss of generality, we may
assume
$x=(0,0)$,
$(0, -r) \in \partial D$ with $B((0, -2r), r)
 \in \bR^2 \setminus  D$.

Let $\wh a := ar^{(2-\alpha)/\alpha}$, $\wh y :=r^{-1} y$ and $\wh D
:=r^{-1} D$. Then by \eqref{e:scaling2},
\begin{equation}\label{e:sc1}
G^{a}_D ( 0, y) = G^{\wh a}_{\wh D} (0, \wh y).
\end{equation}
By the strong Markov property, we have
\begin{equation}\label{e:sc2}
G^{\wh a}_{\wh D}(0, \wh y)=G^{\wh a}_{B(0,1)} (0, \wh y) +\E_0
\left[ G^{\wh a}_{\wh D}\big(X^{\wh a}_{\tau^{ \wh a}_{B(0,1)}}, \wh
y\, \big) \right].
\end{equation}
Note that $\wh a = ar^{(2-\alpha)/\alpha} \le M
R^{(2-\alpha)/\alpha}$. Define
$$
h(z, w):=\E_z \left[ G^{\wh a}_{\wh D}\big (X^{\wh a}_{\tau^{\wh
a}_{B(0,1)}}, w\, \big) \right].
$$
For each fixed $z\in B(0,1)$, the function $w \mapsto
h(z, w)$ is harmonic in $B(0, 1)$ with respect to $X^{\wh a}$ and
for each fixed $w\in B(0,1)$, $z\mapsto h(z, w)$ is  harmonic in
$B(0, 1)$ with respect to process $X^{\wh a}$. So it follows from
Proposition \ref{uhp},
$$
h(0, \wh y)
 \leq c_1 \, \min_{z, w\in B(0,5/6)} h(z, w) \leq  c_1 \,
\min_{z, w\in B(0,5/6)} G^{\wh a}_{\wh D}(z, w) \le c_1 G^{\wh
a}_{\wh D}(0, x_1)
$$
where $|x_1|=1/2$. In the second inequality we used that $G^{\wh
a}_{\wh D}(\cdot, w)$ is superharmonic in $B(0,1)$
for $X^{\wh a}$. Note that $\wh D
\subset E:=\bR^2 \setminus B((0,-2), 1)$. Thus by Lemma
 \ref{l:Cap2},
\begin{equation}\label{e:sc3}
h(0, \wh y) \leq c_1  G^{\wh a}_{\wh D}(0, x_1)  \le c_1 G^{\wh
a}_{E}(0, x_1) \le
\frac{c_2}{\mbox{\rm Cap}^0_{E} \big(\overline{B(0, 1/4)}\big)} <\infty .
\end{equation}
On the other hand, by Corollary  \ref{c:Cap2}
\begin{equation}\label{e:sc4}
G^{\wh a}_{B(0,1)} (0, \wh y) \le c_3 \log \left( \frac{1}{|\wh
y|}\right)= c_3 \log \left( \frac{\delta_D(0)}{|y|}\right).
\end{equation}
It follows from \eqref{e:sc1}-\eqref{e:sc4} that
$$ G_D^a(x, y)=G_D^a(0, y) \leq c_4+c_3\log \left( \frac{\delta_D(0)}{|y|}\right)
\leq c_5 \log \left( \frac{\delta_D(x)}{|x-y|}\right),
$$
which proves the lemma.
\qed

\begin{lemma}\label{l:Cap4}
Suppose that $M>0$ and that $D$ is a bounded $C^{1, 1}$ open set in
$\R^2$. If $x$ and $y$ are in the same component of $D$ with
$$
\frac{1}{c} (\delta_D(x) \vee \delta_D(y)) \le |x-y| \le c
(\delta_D(x) \wedge  \delta_D(y))
$$
for some $c>1$, then there exists $C_{17}=C_{17}(c, D)>0$ such that
for every $a \in [0,M]$
$$
G^a_{D}(x,y)\,\le\, C_{17}.
$$
\end{lemma}
\pf Without loss of generality, we assume $\delta_D(x) \le
\delta_D(y)$. If $\frac{1}{c} \delta_D(y) \le |x-y| \le
\frac{3}{4} \delta_D(x) $, then the lemma follows from Lemma
\ref{l:Cap3}. In the case $\frac{3}{4} \delta_D(x) \le |x-y| \le c
\delta_D(x) $,
since $x, y$ are in the same component of $D$,
we use Proposition \ref{uhp} and a standard Harnack
chain argument. \qed

\begin{thm}\label{t:Gest20}
Suppose that $M>0$ and that $D$ is a bounded $C^{1, 1}$ open set in
$\R^2$. There exists $C_{18}=C_{18}(D)>0$ such that for every $a \in
(0,M]$ and all $x, y \in D$
$$
G^a_D(x,y)\,\le\, C_{18}\,\log \left(1+\frac1{|x-y|^2} \right).
$$
\end{thm}

\pf
Let $L=\max\{2\, \mathrm{diam}(D),2\}$.

\medskip \noindent
(i) If $|x-y| < 1/4$, by Lemma \ref{l:Cap3}
applied to $B(x,L)$,
$$
G^a_D(x,y)\,\le\,G^a_{B(x,L)}(x,y) \le c_1\,\log \left(\frac1{|x-y|}
\right) \le c_2\,\log \left(1+\frac1{|x-y|^2} \right).
$$

\medskip \noindent
(ii) If $1/4 \le|x-y|$, by \eqref{e:scaling2} and
Corollary \ref{c:Cap2},
\begin{eqnarray*}
 G^a_D(x,y) &\le& G^a_{B(x,L)}(x,y)\le G^{a
L^{(2-\alpha)/\alpha}}_{B(L^{-1}x,1)}\big(L^{-1}x,L^{-1}y\big)
\leq c_3\,\log \left( \frac{L}{|x-y|} \right) \\
& \le &  c_3\,\log \left(
4L \right)  \le  c_4 \le  c_5\,\log \left(1+\frac1{|x-y|^2} \right).
\end{eqnarray*}
\qed

Now we are ready to prove Theorem \ref{t-main-green} for  $d=2$.
Recall that $f_D(x, y)$ is defined in \eqref{e:f1}.

\pff {\bf of Theorem \ref{t-main-green} when $d=2$:} By Theorem
\ref{T:gfcesks2}, we only need to consider
the upper bound. We divide its proof into two steps.

{\it Step 1}. We first consider the case that $x$ and $y$ are in the
same component of $D$. Without loss of generality, throughout
 this proof, we assume $\delta_D(x) \le \delta_D(y)$.

Fix $z_0\in D$ with $\kappa R<\delta_D(z_0)<R$, and let
$\eps_1:=\kappa R/24$. Choose $Q_x, Q_y \in \partial D$ with
$|Q_x-x| =\delta_D(x)$ and $|Q_y-y| =\delta_D(y)$. We consider the
following five cases separately.

\medskip \noindent
(a)
If $\delta_D(x ) \ge \eps_1 \kappa^2/32$,  by Theorem \ref{t:Gest20}
$$
G^a_D(x,y)\,\le\,c_1 \log \left(1+\frac1{|x-y|^2} \right) \,
\le\, c_2 \log \left(1+f_D(x,y) \right).
$$

\medskip \noindent
(b)
Suppose $\delta_D(x ) < \eps_1 \kappa^2/32$ and $\delta_D(y ) \ge
\eps_1 \kappa/4$.
Let
$r:= \eps_1 \kappa/16$ and put $x_1=A_{r\kappa /2}(Q_x)$.
One can easily check that $|z_0-Q_x| \ge r$ and $|y-Q_x| \ge r$. So
by (\ref{e:CG_1}), Theorem \ref{t:main} and Theorem \ref{t:Gest20},
we have
$$
G^a_D(x,y) \,\le\, c_3 G^a_D(x_1,y) \frac{G^a_D(x, z_0)}
{G^a_D(x_1, z_0)} \le c_4 \delta_D(x)  \le c_5 f_D(x,y)
$$
for some $c_3, c_4, c_5>0$. Note that  $f_D(x,y)< c_6$ in this case
because $D$ is bounded and $|x-y| \ge \delta_D(y)-\delta_D(x) \ge
\eps_1 \kappa(1/4 -\kappa/32)>0$.
 So it follows from the above display and  Lemma \ref{l:Gest0} that
$$
G^a_D(x,y) \,\le\, c_7 \log \left(1+f_D(x,y) \right).
$$

\medskip \noindent
(c) Suppose $\delta_D(x ) < \eps_1 \kappa^2/32$, $\delta_D(y) \le
\eps_1 \kappa/4$ and $|x-y| < \delta_D(y)/2$. From $\delta_D(y)\le
|x-y|+\delta_D(x)$ we conclude that $\delta_D(y)<2\delta_D(x)$ and
so $|x-y|<\delta_D(x)$. This together with Lemma \ref{l:Cap3}
gives that
$$
G^a_{D}(x,y)\,\le\,
c_8\,\log \left(\frac{\delta_D(y)}{|x-y|} \right)\,\le\,
c_9\,\log \left(1+f_D(x,y) \right).
$$

\medskip \noindent
(d)
If $\frac12 \delta_D(y) \le |x-y| \le(24/\kappa^2)\delta_D(x),$ by Lemma
\ref{l:Cap4},
$$
G^a_D(x,y) \le c_{10} \le c_{11} \log (1+f_D(x,y)).
$$

\noindent (e) The remaining case is
$$
\delta_D(x) \le \frac{\eps_1 \kappa^2}{32}, \quad \delta_D(x) \le
\delta_D(y)\le\frac{\eps_1\kappa}4 \quad \hbox{and} \quad  |x-y| >
\max\left\{\frac{24}{\kappa^2}\delta_D(x) , \frac{\delta_D(y)}{2}\right\}.
$$
We claim that in this case
\begin{equation}\label{e:Gest2}
G^a_D(x,y)\,\le\, c_{12}\,f_D(x,y).
\end{equation}
By Lemma \ref{l:Gest0}, the above implies that $G_D^a(x, y) \leq
c_{13} \log (1+f_D(x, y))$ since in this case
$$ f_D(x, y) = \frac{\delta_D(x)\delta_D(y)}{|x-y|^2} \leq 4.
$$
We now proceed
to prove \eqref{e:Gest2} by considering the following two subcases.

\medskip \noindent
(i) $(24/\kappa^2)\delta_D(x) \le |x-y| \le (4/ \kappa)\delta_D(y)$:
Let $r:= \delta_D(y)/3$.
Put $x_1=A_{r\kappa/2}(Q_x)$. One can easily check that $|z_0-Q_x|
\ge r$ and $|y-Q_x| \ge r$. So by (\ref{e:CG_1}) and
Theorem \ref{t:main}, we have
$$
G^a_D(x,y) \,\le\, c_{14} G^a_D(x_1,y) \frac{G^a_D(x, z_0)}{G^a_D
(x_1, z_0)} \le c_{15} G^a_D(x_1,y) \frac{\delta_D(x) }{|x-y|}.
$$
Moreover,
$$
\frac{24}{\kappa^2}\delta_D(x) \le |x-y| \le \frac32|x_1-y| \le
\frac{6}{\kappa}
\delta_D (y) \le \frac{72}{\kappa^3}\delta_D(x_1),
$$
implying that
$$
\delta_D(y)\le |x-y|+\delta_D(x)\le \left(\frac32 +\frac{36}{\kappa^2}\right)|x_1-y|\, .
$$
It follows from Lemma \ref{l:Cap4} that $G^a_D(x_1,y) \le c_{16}$. Therefore
$$
G^a_D(x,y) \,\le\, c_{17}  \frac{\delta_D(x) }{|x-y|} \,\le\, c_{18}
f_D(x, y).
$$

 \medskip \noindent
(ii) $\delta_D(x) \le \delta_D(y)\le (\kappa/4) |x-y|$:
Let $r:= \frac12 (|x-y| \wedge
\eps_1)$.
Put $x_1=A_{r\kappa/2}(Q_x)$ and $y_1=A_{r\kappa/2}(Q_y)$.
 Then, since $|z_0-Q_x| \ge r$ and $|y-Q_x| \ge r$, by (\ref{e:CG_1}), we
have
$$
c_{19}^{-1}\,  \frac{G^a_D(x_1,y)}{G^a_D(x_1, z_0)}\,\le\,
\frac{G^a_D(x,y)}{G^a_D(x, z_0)} \,\le\, c_{19}
\frac{G^a_D(x_1,y)}{G^a_D(x_1, z_0)}
$$
for some $c_{19}>1$. On the other hand, since  $|z_0-Q_y| \ge r$ and
$|x_1-Q_y| \ge r$, applying  (\ref{e:CG_1}),
$$
c_{19}^{-1}\,\frac{G^a_D(x_1,y_1)}{G^a_D(x_1, y)} \,\le\,
\frac{G^a_D(y_1,z_0)}{G^a_D(y, z_0)} \,\le\, c_{19}
\frac{G^a_D(x_1,y_1)}{G^a_D(x_1, y)}.
$$
Putting the four inequalities above together we get
$$
c_{19}^{-2}\,\frac{G^a_D(x_1,y_1)}{G^a_D(x_1, z_0)G^a_D(y_1, z_0)}
\,\le\, \frac{G^a_D(x,y)}{G^a_D(x, z_0)G^a_D(y, z_0)} \,\le\,
c_{19}^2 \, \frac{G^a_D(x_1,y_1)}{G^a_D(x_1, z_0)G^a_D(y_1, z_0)}.
$$
Moreover, $\frac13|x-y| < |x_1-y_1| < 2 |x-y| $ and
$$
\frac{4}{3\kappa}(\delta_D(x_1) \vee \delta_D(y_1))  \le \frac13
|x-y| \le |x_1 -y_1| \le \frac{64}{\kappa^3\eps_1} (\delta_D(x_1) \wedge
\delta_D(y_1)).
$$
Thus by Lemma \ref{l:Cap4} and Theorem \ref{t:main}, we have
$$
G^a_D(x,y) \,\le\,c_{20} \frac{G^a_D(x, z_0)G^a_D(y,
z_0)}{G^a_D(x_1, z_0) G^a_D(y_1, z_0)}\, \le c_{21} f_D(x,y)
$$
for some $c_{20}, c_{21}>0$. This completes the proof of the claim
\eqref{e:Gest2} and therefore of the theorem when $x$ and $y$ are
are in the same component of $D$.

\medskip

{\it Step 2.} Next we consider the case that $x$ and $y$ are in two
different components of $D$. This part of the proof is the same as
the second part of the proof of Theorem \ref{t-main-green} when
$d\geq 3$ (that is, the paragraph containing
\eqref{e:3.12}-\eqref{e:3.13}). The only place that needs
modification is the proof of $\sup_{z\in D} \E_x [\tau^a_D]\leq
c_{22}<\infty$. When $d=2$, we can not use \eqref{ub} to deduce it.
However, since $D$ is bounded, there is $K>0$ so that $D\subset B(0,
K)$. It follows from {\it Step 1} that
$$
\sup_{z\in D, a\in(0, M]} \E_z [\tau^a_D] \leq \sup_{z\in B(0, K),
a\in(0, M]} \E_z [ \tau^a_{B(0, K)}]
\leq c_{23}<\infty.
$$
This completes the proof of Theorem \ref{t-main-green}. \qed

\medskip

\begin{thm}[3G theorem for $d=2$]\label{t:3Gt2}
Suppose that $M>0$ and that $D$ is a bounded  $C^{1,1}$ open set in
$\R^2$. Then there exist positive constants $C_{19}=C_{19}(D,
\alpha, M)$ and $C_{20}=C_{20}(D, \alpha, M)$ such that for all $x,
y, z \in D$ and $a \in (0, M]$
\begin{eqnarray*}
\frac{G^a_D(x,y) G^a_D(y,z)} { G^a_D(x,z)} & \le&
C_{19}\, \left(\log(1+f_D(x,y)) + \log(1+f_D(y,z)) +1 \right) \\
&\le& C_{20} \left( \big(\log\frac1{|x-y|} \vee 1 \big) +
\Big(\log\frac1{|y-z|} \vee 1 \Big)\right).
\end{eqnarray*}
\end{thm}

\pf Note that, if $x,z$ are in different components of $D$, either
$x,y$ or $y,z$ are in different components of $D$. Thus, by Theorem
\ref{t-main-green} for $d=2$ and the fact that $a \in [0, M]$, we
have
$$
\frac{G^a_D(x,y) G^a_D(y,z)} { G^a_D(x,z)} \le
c_1 \frac{
\log(1+f_D(x,y))\log(1+f_D(y,z))}{
\log(1+f_D(x,z))}
$$
for some $c_1=c_1(M, D, \alpha)$.
Now following the proof of \cite[Theorem 6.24]{CZ}, we get the
theorem.
\qed

\begin{remark} \rm
By considering  how many different components of $D$ that $x, y$ and
$z$ fall into, we could get more precise 3G estimates with the
dependence on $a$ explicitly spelled out. Theorem \ref{t:3Gt2} will
not be used in the remainder of this paper.
\end{remark}

\section{One Dimensional Case}

In this section we assume $d=1$ and prove Theorem \ref{t-main-green}
for this case. We will follow the ideas in \cite{KSV09}.

Let ${\overline X}^a$ be the supremum process of $X^a$ defined by
${\overline X}^a_t=\sup\{0\vee {\overline X}^a_s:0\le s\le t\}$ and
let ${\overline X}^a-X^a$ be the reflected process at the supremum.
The local time at zero of ${\overline X}^a-X^a$ is denoted by
$L^a=(L^a_t: t\ge 0)$ and the inverse local time by $\{\tau^a_t:
t\geq 0\}$, where $\tau^a_t:=\inf\{s: L^a_s > t\}$. The inverse
local time $\{\tau^a_t: t\geq 0\}$ is a subordinator. The
(ascending) ladder height process of $X^a$ is the process
$H^a=(H^a_t: t\ge 0)$ defined by $H^a_t=X^a_{\tau^a_t}$. The ladder
height process is again a subordinator. It follows from \cite{KSV09}
that $H^a$ is a special subordinator with Laplace exponent given by
\begin{equation}\label{e:Fri}
\chi^a(\lambda)=\exp\left(\frac1\pi\int^\infty_0\frac{
\log(\theta^2\lambda^2+ a^{\alpha}\theta^\alpha\lambda^\alpha)}
{1+\theta^2}d\theta\right)
\end{equation}
and that the drift coefficient of $H^a$ is 1. When $a=0$, we have
$\chi^0(\lambda)=\lambda$. Thus, if $V^a$ is the potential measure
of $H^a$ and $V^a(x)=V^a([0, x])$, then, for every $a\geq 0$, $V^a$
has a continuous, decreasing and strictly positive potential density
$v^a$ such that $v^a(0+)=1$. When $a=0$, we have $v^a\equiv 1$. The
following results is a uniform version of \cite[Proposition
2.3]{KSV09} in our present special case.

\begin{lemma}\label{l:estimate-for-V}
Let $M$ and $R_2$ be positive constants. There exists a constant
$C_{21}=C_{21}(M, R_2)\in (0,1)$ such that for all $a\in [0, M]$ and $x\in
(0,R_2]$,
$$
C_{21} \le v^a(x) \le C_{21}^{-1} \quad \text{and}\quad C_{21} x\le V^a(x)
\le C_{21}^{-1}x\, .
$$
\end{lemma}

\pf Since $H^a$ is special, the potential density $v^a$ is a
decreasing function. Hence $\inf_{0<t\le R_2} v^a(t)=v^a(R_2)$. It
follows from
\eqref{l:estimate-for-V} that the Laplace exponent
$\chi^a$ is continuous in $a$. Thus, the potential measures converge
vaguely, and by continuity and monotonicity of $v^a$, we get that
$v^a(t)\to v^b(t)$ as $a\to b$ for all $t>0$. In particular,
$v^a(R_2)\to v^b(R_2)$. Therefore $c_1:=\inf_{0<t\le R_2, 0\le a \le
M}v^a(t) >0. $ Since $v^a(t)\le 1$ for all $t>0$ and all
$a \ge 0$,
we
get that $c_2 =\sup_{0<t\le R_2, 0\le a \le M}v^a(t) =1$. Choose
$c_3=c_3(M, R_2)\in (0,1)$ such that $c_3\le c_1\le c_2 \le
c_3^{-1}$. Since $V(x)=\int_0^x v(t)\, dt$, the claim follows
immediately. \qed

\begin{thm}\label{T:gfcesks1d}
Suppose that $M>0$. For any bounded open interval $D$ in $\R$, there
exists a constant $C_{22}=C_{22}(\alpha, M, D)
>1$ such that for all
$x,y \in D$ and all $a\in (0, M]$,
$$
C_{22}^{-1}
\left(
\left(\delta_D(x)\delta_D(y)\right)^{1/2}
\wedge\frac{\delta_D(x)\delta_D(y)}{ |x-y|}
\right)
\le   G^a_D(x,y) \le
C_{22}\,
\left(
\left(\delta_D(x)\delta_D(y)\right)^{1/2}
\wedge\frac{\delta_D(x)\delta_D(y)}{ |x-y|}
\right).
$$
\end{thm}

\pf
The proof of the lower bound is similar to that of Theorem
\ref{T:gfcesks} and \cite[Proposition 3.3]{KSV09}. Using our Lemma
\ref{e:Fri} instead of \cite[Proposition 2.3]{KSV09}, we can follow the
proof of \cite[Propsoition 3.1]{KSV09} to get the upper bound. We
omit the details. \qed

Integrating the estimate in Theorem \ref
{T:gfcesks1d} with respect to
$y$ yields the following corollary.

\begin{cor}\label{c-lb-exit1d}
Suppose that $M>0$. For any bounded open interval $D$ in $\R$, there
exists a positive constant
$C_{23}=C_{23}(\alpha, M, D)>1$ such that
for all $x\in D$ and all $a\in (0, M]$,
$$
C_{23}^{-1}\delta_D(x)\le \E_x [\tau_D^a]\le C_{23}
\delta_D(x)\, .
$$
\end{cor}

Using Corollary \ref{c-lb-exit1d},
we can repeat the proof of
Theorem \ref{t-main-green} for $d\geq 3$  case (see also
\cite[Theorem 3.8]{KSV09}) to generalize
Theorem \ref{T:gfcesks1d}
to general (not necessarily connected) bounded $C^{1,1}$ open sets.
Recall that $g_D^a$ is defined by \eqref{e:1.21d}.

\begin{thm}\label{T:gfcesks21d}
Suppose that $D$ is a bounded $C^{1, 1}$ open set in $\bR$. There
exists a positive constant
$C_{24}=C_{24}(\alpha, M, D)>1$
such that
for all $x,y \in D$ and all $a\in (0, M]$,
$$
C_{24}^{-1}\, g_D^a (x, y) \le G^a_D(x,y) \le  C_{24}
\, g_D^a (x,
y).
$$
\end{thm}

\section{Martin Boundary and Martin Kernel Estimates}

Throughout this section we assume that $d\geq 1$ and $D$ is a
bounded $C^{1,1}$ open set in $\RR^d$ with characteristics $(R,
\Lambda)$ and the corresponding $\kappa$. We will show in this
section that the Martin boundary and  the minimal Martin boundary of
$D$ with respect to $X^a$ can both be identified with the Euclidean
boundary $\partial D$ of $D$. With the boundary Harnack principle
given in Theorem \ref{t:main}, the arguments of this section are
modifications of the corresponding parts of
\cite{Bo, CS, KS2, KSV}.
For this reason, most of the proofs in this section will be omitted.

The next lemma follows from Theorem \ref{t:main}.

\begin{lemma}\label{l:5B}
Suppose that $M>0$ and that $D$ is a bounded $C^{1, 1}$ open set in
$\R^d$. There exists a positive constant $C_{25}=C_{25}(D, \alpha,
M)$ such that for all $a \in (0,M]$, $Q\in \partial D$, $r\in (0,
R/2)$, and nonnegative function $u$ in $\R^d$ which is
 harmonic with respect to $X^a$ in $D \cap B(Q, r)$ we have
\begin{equation}\label{e:gamma}
u(A_r(Q))\,\le\, C_{25}\,\left(\frac2{\kappa}\right)^{k}\,
u(A_{(\kappa/2)^{k}r}(Q)), \qquad k=0, 1, \dots .
\end{equation}
\end{lemma}

\begin{lemma}\label{L:2.00}
Suppose that $M>0$. For every $b \in (0,\infty)$, there exist
$C_{26}=C_{26}(M, b)>0$ and $C_{27}=C_{27}(M, b)>0$ such that for
all $x_0 \in \bR^d$, $a \in (0,M]$ and $r \in (0, b]$,
 \bee \label{e:ext}
C_{26} r^2 \, \le \, \E_{x_0}\left[\tau^a_{B(x_0,r)}\right]\, \le\,
C_{27}\,r^2
 \eee
 and
 \bee \label{e:ext2}
\E_{x}\left[\tau^a_{B(x_0,r)}\right]\, \le\,
C_{27} \,r
\delta_{B(x_0,r)}(x).
 \eee
\end{lemma}
\pf See \cite[Lemmas 2.2 and 2.3]{SV05} or \cite[Lemmas 2.3 and
2.4]{CK08} for a proof of \eqref{e:ext}. The inequality
\eqref{e:ext2} follows easily from Theorem \ref{t-main-green}. In
fact, by \eqref{e:scaling1} and Theorem \ref{t-main-green} (with
$Mb^{(2-\alpha)/\alpha}$ instead of $M$)
$$
\E_{x}\left[\tau^a_{B(0,r)}\right]=r^{2}\int_{B(0,1)}
G^{ar^{(2-\alpha)/\alpha}}_{B(0,1)}
(r^{-1} x, z)dz \, \le \, c \, r^{2} \delta_{B(0,1)}(r^{-1} x )
\, =\, c\,  r\,  \delta_{B(0,r)}(x).$$
\qed

For an open set $U\subset \R^d$, let
\begin{equation}\label{PK}
K^a_U(x,z)\,:=\,
\int_U
{G^a_U(x,y)}J^a(y,z) dy, \quad (x,z) \in U \times
\overline{U}^c.
\end{equation}
Then by \eqref{e:levy}, for any non-negative measurable function $f$
on $\bR^d$,
$$
\E_x\left[f(X^a_{\tau^a_U});\,X^a_{\tau^a_U-} \not= X^a_{\tau^a_U}
\right] =\int_{\overline{U}^c} K^a_U(x,z)f(z)dz.
$$

From \eqref{PK}, Theorem \ref{t-main-green}  and Lemma \ref{L:2.00},
we immediately get the following proposition.

\begin{prop}\label{p:Poisson1}
Suppose that $M>0$. There exist $C_{28}>0$ and $C_{29}>0$ such that
for all $a \in (0,M]$ and $r \in (0, R)$ and $x_0 \in \R^d$,
\begin{equation}\label{e:P1}
K^a_{B(x_0,r)}
(x,y) \,\le\, C_{28} \,r\,
(r-|x-x_0|)(|y-x_0|-r)^{-d-\alpha} \qquad \hbox{for }  (x,y) \in
B(x_0,r)\times \overline{B(x_0,r)}^c
\end{equation}
and
\begin{equation}\label{e:P2}
K^a_{B(x_0,r)}
(x_0,y) \,\ge\, C_{29}\, r^2\, |y-x_0|^{-d-\alpha}
\qquad \hbox{for } y \in \overline{B(x_0,r)}^c.
\end{equation}
\end{prop}

Using \eqref{e:P2}, the proof of the next lemma is similar to that
of \cite[Lemma 4.3]{KS2} or \cite[Lemma 5.3]{KSV}. Thus we skip the
proof.

\begin{lemma}\label{l:la}
Suppose that $a >0$ and that $D$ is a bounded $C^{1, 1}$ open set in
$\R^d$. There
exists a positive constant $C_{30}=C_{30}(D, \alpha, a)$
such that for all $Q \in \partial D$, $r \in (0, R/2)$ and $w \in
D\setminus B(Q, r)$,
$$
G^a_D(A_r(Q), w)\, \ge\, C_{30}\, r^{2} \int_{B(Q, r)^c}
 j^a(|z- Q|/2) G^a_D(z, w)dz.
$$
\end{lemma}

Using \eqref{e:P1}, Lemmas \ref{l:5B} and \ref{l:la}, the proof of
the next lemma is similar to that  of \cite[Lemma 4.4]{KS2} or
\cite[Lemma 5.4]{KSV}. Thus we skip the proof.

\begin{lemma}\label{l:14B}
Suppose that $a >0$  and that $D$ is a bounded $C^{1, 1}$ open set
in $\R^d$. There exist positive constants
$C_{31}=C_{31}(D,\alpha,a)$ and $C_{32}=C_{32}(D,\alpha,a)<1$ such
that for any $Q \in \partial D $, $r\in (0, R/4)$ and  $w \in D
\setminus B(Q,2r/\kappa)$,  we have
$$
\E_x\left[G^a_D(X^a_{\tau^a_{D \cap B_k}}, w):\,X^a_{\tau^a_{D \cap
B_k}} \in B(Q, r)^c  \right] \,\le\, C_{31}\,C_{32}^{k} \,
G^a_D(x,w), \quad x \in D \cap B_k,
$$
where $B_k:=B(Q, (\kappa/2)^{k}r)$, $ k=0,1, \dots$.
\end{lemma}

Let $x_0\in D$ be fixed and set
$$
M^a_D(x, y):=\frac{G^a_D(x, y)}{G^a_D(x_0, y)}, \qquad x, y\in D,~
y\neq x_0.
$$

Now the next theorem follows from Theorem \ref{t:main} and Lemma
\ref{l:14B}
  (instead of  \cite[Lemma 13]{Bo} and \cite[Lemma 14]{Bo}, respectively)
  in very much the same way as in the case of symmetric
stable processes in
\cite[Lemma 16]{Bo} (with Green functions
instead of harmonic functions). We omit the details.

\begin{thm}\label{t2.2}
Suppose that $a >0$  and that $D$ is a bounded $C^{1, 1}$ open set
in $\R^d$. There exist positive constants $R_1$, $M_1$, $C_{33}$ and
$\beta$ depending on $D$, $\alpha$ and $a$ such that for any $Q \in
\partial D $, $r < R_1$ and $z \in D \setminus B(Q, M_1 r)$, we have
$$
\left|M^a_D(z, x)-M^a_D(z, y)\right| \,\le\,
C_{33}\,\left(\frac{|x-y|}r\right)^{\beta}, \qquad x, y\in  D \cap
B(Q, r).
$$
In particular, the limit $\lim_{D \ni y\to w} M^a_D(x, y)  $ exists
for every $w\in \partial  D$.
\end{thm}

As the process $X^{a,D}$ satisfies Hypothesis (B) in Kunita and
Watanabe \cite{KW}, the process $X^{a,D}$ has a Martin boundary: For
every $a \in (0,M]$,  there is a compactification $D_a^M$ of $D$,
unique up to a homeomorphism, such that $M^a_D(x, y)$ has a
continuous extension to $D\times (D_a^M\setminus\{x_0\})$ and
$M^a_D(\cdot, z_1) =M^a_D(\cdot, z_2)$ if and only if $z_1=z_2$
(see, for instance,  \cite{KW}). The set
$\partial_a^MD=D_a^M\setminus D$ is called the Martin boundary of
$D$ for $X^{a,D}$. For $z\in
\partial_a^MD$, set $M^a_D(\cdot, z)$ to be zero in $D^c$.

For each fixed $z\in \partial D$ and $x\in D$, let
$$
M^a_D(x, z):=\lim_{D\ni y\to z}M^a_D(x,y),
$$
which exists by Theorem \ref{t2.2}.
$M^a_D$ is called the Martin kernel of $D$ with respect to $X^a$.
For each $z\in \partial D$, set
$M^a_D(x, z)$ to be zero for $x\in D^c$. By Theorem \ref{t2.2},
$M^a_D(z, x)$ is  jointly continuous on $\{x\in D: \delta_D(x) >
2\eps \} \times \partial D$,  and hence on $D\times \partial D$
after letting $\eps \downarrow 0$.

The following Martin kernel estimate is an immediate consequence of
the Green function estimates in Theorem \ref{t-main-green}. Recall
that $D(x)$ denotes the  component of $D$ that contains $x$. Define
\begin{equation}\label{e:5.7}
h_D^a(x, z):= \begin{cases} \frac {\delta_D (x)} {|x - z |^{d}}
\quad &\hbox{if } x\in D(x_0), z\in \partial D(x_0)
\hbox{ or } x\in D\setminus D(x_0), z\in \partial  D \setminus
(\partial D(x_0)\cup \partial D(x)), \\
a^\alpha  \,\frac {\delta_D (x)} {|x - z |^{d}}
\quad &\hbox{if } x\in D(x_0), z\in \partial D
\setminus \partial D(x_0),\\
a^{-\alpha} \, \frac {\delta_D (x)} {|x - z |^{d}}
\quad &\hbox{if } x\in D\setminus D(x_0), z\in \partial D(x_0).
\end{cases}
\end{equation}

\begin{thm}\label{T:4.2}
Suppose that $M>0$ and that $D$ is a bounded $C^{1, 1}$ open set in
$\R^d$. There exists $C_{34}:=C_{34}(x_0, D, \alpha,M)
>1$ such that
for all $a \in (0, M]$,
$$
C_{34}^{-1} \, h_D^a(x, z) \leq M^a_D (x, z) \leq C_{34} \, h_D^a(x,
z) \qquad \hbox{for } x\in D, z\in \partial D.
$$
\end{thm}

\medskip

Theorem \ref{T:4.2} in particular implies that $M^a_D(\cdot , z_1)$
differs from $M^a_D(\cdot, z_2)$ if $z_1$ and $z_2$ are two
different  points on $\partial D$.

Now using our Green function estimates, \eqref{e:P1} and Lemma
\ref{l:5B}, one can follow the arguments in the proofs of
\cite[Lemmas 4.6--4.7]{KS2} and \cite[Lemmas 5.6--5.7]{KSV} and get
the next two lemmas.

\begin{lemma}\label{l:MH1}
Suppose that $D$ is a bounded $C^{1, 1}$ open set in $\R^d$.
For every $z \in \partial D$ and $B \subset \overline{B} \subset D$,
$M^a_D(X^a_{\tau^a_{B}} , z)$ is $\P_x$-integrable.
\end{lemma}

\begin{lemma}\label{l:MH2}
Suppose that $D$ is a bounded $C^{1, 1}$ open set in $\R^d$.
For every $z \in \partial D$ and  $x \in D$,
\begin{equation}\label{eqn:4.1}
M^a_D(x,z) \,=\, \E_x \left[M^a_D\big(X^D_{\tau^a_{B(x,r)}},
z\big)\right]
 \quad \mbox{ for every } 0<r\le \frac12(R\wedge
\delta_D(x)).
\end{equation}
\end{lemma}

Unlike the case in the proofs of
\cite[Theorem 2.2]{CS} and
\cite[Theorem 4.8]{KS2},
$\P_x(X^a_{\tau^a_{U}} \in \partial
U)\not=0$ for every smooth open set $U$. Thus we give the details of
the proof of the next result.

\begin{thm}\label{T:L4.3}
For every $z \in \partial D$, the function $x\mapsto M^a_D(x, z)$ is
harmonic in $D$ with respect to $X^a$.
\end{thm}

\pf Fix $z\in \partial D$ and let $h(x):=M_D^a(x,z)$. Consider an
open set $D_1 \subset\overline{D_1}\subset D$ and $x \in D_1$ and
put
$$
r(x)= \frac12(R\wedge \delta_D(x)) \quad \mbox{ and }\quad B(x)=B(x,
r(x)).
$$
Define a sequence of stopping times $\{T_m, m\ge1\}$ as follows:
$$
T_1=\inf\{t>0: X^a_t\notin B(X^a_0)\},
$$
and for $m\ge 2$,
$$
T_m= \begin{cases}  T_{m-1}+
\tau_{B(X^a_{T_{m-1}})}\circ\theta_{T_{m-1}}
\qquad &\hbox{if }X^a_{T_{m-1}}\in D_1\\
\tau^a_{D_1} \qquad &\hbox{otherwise}.
\end{cases}
$$
Note  that $X^a_{\tau^a_{D_1}}\in \partial D_1$ on
$\cap_{n=1}^{\infty}\{T_n < \tau^a_{D_1}\}. $ Thus, since $\lim_{m
\to  \infty} T_m=\tau^a_{D_1}$ $\P_x$-a.s. and $h$ is continuous in
$D$,
$$
\lim_{m \to  \infty} h(X^a_{T_m})= h(X^a_{\tau^a_{D_1}}), \qquad
\hbox{on } \cap_{n=1}^{\infty}\{T_n < \tau^a_{D_1}\}
$$
and, since $h$ is bounded on $\overline{D_1}$,  by the dominated
convergence theorem
\begin{eqnarray*}
\lim_{m\to\infty} \E_x\left[h(X^a_{T_m});\, \cap_{n=0}^{\infty}
\{T_n<\tau^a_{D_1} \}\right]= \E_x\left[h(X^a_{\tau^a_{D_1}});\,
\cap_{n=0}^{\infty} \{T_n<\tau^a_{D_1} \}\right].
\end{eqnarray*}

Therefore, using Lemma \ref{l:MH2}
\begin{eqnarray*}
h(x)& =& \lim_{m\to\infty}\E_x\left[h(X^a_{T_m})\right] \\
&=& \lim_{m\to\infty}\E_x\left[h(X^a_{T_m});\, \cup_{n=0}^{\infty}
\{T_n=\tau^a_{D_1} \}\right] + \lim_{m\to\infty} \E_x
\left[h(X^a_{T_m});\, \cap_{n=0}^{\infty}
\{T_n<\tau^a_{D_1} \}\right]\\
&=& \E_x\left[h(X^a_{\tau^a_{D_1}});\,
\cup_{n=0}^{\infty}
\{T_n=\tau^a_{D_1} \}\right] +  \E_x\left[h(X^a_{\tau^a_{D_1}});\,
\cap_{n=0}^{\infty}
\{T_n<\tau^a_{D_1} \}\right]\\
& = & \E_x\left[h(X^a_{\tau^a_{D_1}})\right].
\end{eqnarray*}
\qed

\medskip

A consequence of Theorems \ref{t2.2}, \ref{T:4.2} and \ref{T:L4.3}
is that, when $D$ is a bounded $C^{1, 1}$ open set, the Martin
boundary of $X^{a,D}$ can be identified with the Euclidean boundary
$\partial D$ of $D$.

A positive harmonic function $u$ for  $X^{a,D}$ is minimal if, whenever
$v$ is a positive harmonic function for $X^{a,D}$ with $v\le u$ on $D$,
one must have $u=cv$ for some constant $c$. The set of points $z\in
\partial_a^MD$ such that $M^a_D(\cdot, z)$ is minimal harmonic
for $X^{a,D}$ is called the minimal Martin boundary of $D$ for $X^{a,D}$.

With the explicit estimates from Theorem \ref{T:4.2}, by the same
argument as that for \cite[Theorem 3.7]{CS}, we have the following.

\begin{thm}\label{T:4.6}
Suppose that $D$ is a bounded $C^{1,1}$ open set in $\R^d$ and
$a>0$.
 For every $z\in \partial D$, $M^a_D(\cdot, z)$ is a minimal harmonic
function for $X^{a,D}$. Thus the minimal Martin boundary of $D$ can
be identified with the Euclidean boundary.
\end{thm}

We know from the general theory in Kunita and Watanabe \cite{KW}
that non-negative superharmonic functions with respect to $X^{a,D}$
(or equivalently, superharmonic functions with respect to $X^{a}$
that vanish on $D^c$) admit a Martin representation. Thus,  by
Theorem \ref{T:4.6} we conclude that,  for every  superharmonic
function $u\geq 0$ with respect to $X^{a,D}$, there is a  unique
Radon measure $\mu$ in $D$ and a finite measure $\nu$ on $\partial
D$ such that
\begin{equation}\label{eqn:martin}
 u(x) = \int_D G^a_D(x, y)\mu (dy)+\int_{\partial D} M^a_D(x,z) \nu (dz).
\end{equation}
Furthermore, $u$ is harmonic for $X^{a,D}$ if and only if the measure $\mu=0$.

\section{Perturbation Results}\label{s:pr}

In this section we assume $d \ge 1$ and fix $a>0$. We consider a
symmetric L\'evy process $Z$ which can be thought of as a
perturbation of $X^a$, and show that under certain conditions, the
Green function of $Z^D$, the process $Z$ killed upon exiting a
bounded $C^{1,1}$ open set $D$, is comparable to the Green function
of $X^{a,D}$, see Theorem \ref{t:gee}. Together with Theorem
\ref{t-main-green}, this gives sharp bounds for the Green function
$G_D^Z$ of $Z^D$.

The approach of this section is motivated by \cite{GR}, where
perturbations of pure jump L\'evy processes are discussed. Even
though they consider pure jump L\'evy processes, some results work
for our case as well.

Throughout this section, $Z$ is a symmetric L\'evy process in
$\bR^d$ such that its L\'evy measure has a density
 $J^Z(x,y)=j^Z(y-x)$ with
respect to the Lebesgue measure. We assume that
$$
j^a_1 (x):= j^a(x)-j^Z(x)
$$
is nonnegative and integrable in $\bR^d$ and put ${\cal
J}^a:=\int_{\R^d} j^a_1(y)dy$. We also assume that the transition
density of
the L\'evy process $Z$
exists and we denote it by $p^Z(t,x,y)=p^Z(t,y-x)$.

Recall that $p^a(t,x,y)=p^a(t,y-x)$ is the transition density
function of $X^a$. It is well known that (see \cite{CK08, SV07})
\begin{equation}\label{e:CK}
p^a(t,x,y) \,\le\, c\, \left(t^{-d/\alpha} \wedge t^{-d/2}  \right)\wedge \left(
t^{-d/2} e^{-c_2|x-y|^2/t}+ \frac{t}{|x-y|^{d+\alpha}}\right),
\quad (t,x,y) \in (0, \infty)
\times \R^d \times \R^d
\end{equation}
for some $c=c(a,d, \alpha)>0$. Thus,  by  following  the  proof of
\cite[Lemma 2.6]{GR}, we have the following.
\begin{lemma}\label{BoundDensity}
$p^Z(t,x)$ is bounded on the set $\{(t,x) : t>0, |x|
> \eps\}$ for $\eps >0$.
\end{lemma}

Recall that
for every bounded open set $D$,
$G^a_{D}$
is the Green function of $X^{a,D}$.
We know from \cite[Corollary 3.12]{KS6}
that there is a constant $ c=c(D,a)$ such that
\begin{equation}\label{e:Ga}
c \, \E_x\big[\tau^a_D\big] \,\E_y\big[\tau^a_D\big]\le  G^a_D(x,y), \quad x,y\in D.
\end{equation}

The proofs of the following three results are the same as those of
\cite[Lemmas 2.2, 2.4, 2.5]{GR}. So we omit their proofs here. In
the remainder of this section, the dependence of the constants on
$Z$ will not be mentioned explicitly.
\begin{lemma}\label{gestzabity2:this}
For every bounded open set $D$, there exists $C_{35}=C_{35}(D, a)>0$
such that for all $x\in D$ and $t\ge 1$ we have
$$
p^a_D(t,x,y) \ \leq \ C_{35} \ t^{-2}  \ \E_x\big[\tau^a_D\big] \
\E_y\big[\tau^a_D\big].
$$
\end{lemma}

For any open set $U\subset \bR^d$, $\tau^Z_U:=\inf\{t>0: \,Z_t\notin
U\}$ denotes the first exit time from $U$ by $Z$.
We
denote by  $Z^{D}$ the subprocess of $Z$ killed upon leaving
$D$ and $p^Z_D(t, x, y)$ the transition density for $Z^D$.

\begin{lemma}\label{momentyogr:this}
For every bounded open set $D$,
there exist a constant $C_{36}=C_{36}(D, a)>1$ such that for every
$x \in D,$
$$
C_{36}^{-1} \ \E_x\big[\tau^Z_D\big] \  \le  \ \E_x \big[\tau^a_D\big] \ \le \ C_{36} \
\E_x\big[\tau^Z_D\big].
$$
\end{lemma}

\begin{lemma}\label{gestzabity1:this}
For every bounded open set $D$ and any $x\in D$ and $t>0$, we have
$$
p^Z_D(t,x,\cdot) \leq  e^{{\cal J}_at} p^a_D(t,x,\cdot)\quad \hbox{a.e.}
$$
If, in addition, we assume that $p^Z(t,\cdot)$ is continuous then we
have for $x,y\in D$,
\bee\label{e:new}
\E_x\left[p^Z(t-\tau^Z_D,Z_{\tau^Z_D},y):t\geq \tau^Z_D \right]\, \leq \, e^{2{\cal
J}_at}\, \E_x\left[p^a(t-\tau^a_D,X^a_{\tau^a_D},y):t\geq \tau^a_D \right].
\eee
\end{lemma}

Using the above lemmas (for Lemma \ref{gestzabity1:this}, only its
first part is needed), and following the proof of \cite[Theorem
3.1]{GR}, we have

\begin{thm}
For every bounded open set $D$,
there exists $C_{37}=C_{37}(D, a)>0 $
 such that
 for every $x \in D$,
\begin{equation}\label{e:G_ub}
G^Z_{D}(x, y)\,\le\, C_{37}\, G^a_{D}(x, y)
\quad \mbox{
a.e. } y \in D.
\end{equation}
\end{thm}

Here are some assumptions that we might put on the process $Z$.

\begin{description}
\item{(A1)}  The
transition density $p^Z_D(t,x,y)$ of $Z^D$ is continuous and
strictly positive in $(0, \infty)\times D \times D$.
\item{(A2)} There exist positive constants $c$ and $\rho$ such that
$j^a_1(x)\le c|x|^{\rho-d}$ on $B(0, 1)$.
\end{description}

\begin{description}
\item{(A3)}  $j^Z$ satisfies either (a) or (b)
below:
\begin{description}
\item{(a)} There exists a non-negative Borel
function $L(x)$, which is  locally integrable on
$\R^d\setminus \{0\}$,
such that for any Borel set $B$,
$
|B|= \int_{B} L(x) j^Z(x) dx.
$
\item{(b)}
There exists $R_0 >0$ such that
$
\inf_{x \in B(0,R_0)}j^Z(x) >0.
$
\end{description}
\end{description}

Without loss of generality, we may and do assume that the constant
$\rho$ is less than 1.

\begin{prop}
[{\cite[Corollary 3.11]{KS6}}]
\label{p:lbG}
Suppose that (A1) and (A3) hold. Then
for every bounded open set $D$,
there exists constant
$C_{38}=C_{38}(D, \alpha)>0$ such that
\begin{equation}\label{lbG}
C_{38}\, \E_x [ \tau^Z_D] \, \E_y [\tau^Z_D] \le\,
G^Z_D(x,y), \qquad \text{for all }(x,y) \in D \times D.
\end{equation}
\end{prop}

For the remainder part of this section,
we assume $D$ is a bounded $C^{1,1}$ open   set
in $\R^d$.

\begin{lemma}\label{G:g2}
Suppose that (A1) and (A3) hold. Then for every $\eps>0$, there
exists $C_{39}=C_{39}(\eps, D,a)>0$ such that for all $x, y \in D$
satisfying $|x-y| \ge \eps$,
$$
G^{a}_D(x,y) \,\le\, C_{39}\,G^Z_D(x,y).
$$
\end{lemma}

\pf
It follows from Theorem \ref{t-main-green} that there exists
$c_1=c_1(D, a)> 0$ such that
\begin{equation}\label{e:g01}
c_1^{-1} \E_{x}\left[\tau^{a}_D\right]   \,\le\,\delta_D(x)\,\le\, c_1
\E_{x}\left[\tau^{a}_D\right] , \qquad x\in D.
\end{equation}
Combining \eqref{e:g01} with Theorem \ref{t-main-green} yields that
there exist $c_2, c_3>0$ so that for all $x$ and $y \in D$ with
$|x-y| \ge \eps$,
$$
G^{a}_D(x,y) \,\le\, c_2\, \delta_D (x) \delta_D (y) \,\le\,
c_3\,\E_{x}[\tau^{a}_D] \E_{y}[\tau^{a}_D].
$$
Therefore, by Lemma
\ref{momentyogr:this}
and Proposition
\ref{p:lbG} we get
$$
G^{a}_D(x,y) \,\le\, c_4 \E_{x}[\tau^Z_D]
\E_{y}[\tau^Z_D]   \,\le\, c_5\,G^Z_D(x,y)
$$
for some positive constants $c_4, c_5>0$. \qed

The next lemma can be proved by following the arguments in the
proofs of \cite[Lemmas 7, 9]{R}. So we skip the details here.

\begin{lemma}\label{l:g5}
Suppose that (A1) holds. For all $x, w\in D$, we have
$$
G^a_D(x,w) \, \le \,  G^Z_D(x,w) \,+\,
\int_D \int_{D}  G^a_D(x,y) j^a_1(y-z)G^a_D(z,w)
dydz.
$$
\end{lemma}

\begin{prop}\label{l:GR2}
Suppose that $d \ge 3$ and that (A2) holds.
 There exists a positive constant $C_{40}=C_{40}(D,
a)$ such that for every $(x,w) \in D \times D$
$$
\int_D\int_D  G^{a}_D(x,y)j^a_1(y-z)  G^{a}_D(z,w) dydz \,\le
\,C_{40}\,G^{a}_D(x,w) |x-w|^{d-2}.
$$
\end{prop}
\pf Using the generalized 3G inequality
(Theorems \ref{t:3Gt})
 for the Green function of $X^{a}$,  one can easily get the following
\begin{eqnarray*}
&& \int_D\int_D  G^{a}_D(x,y)j^a_1(y-z)  G^{a}_D(z,w) dydz \\
&\le &c\,G^{a}_D(x,w) \Big( |x-w|^{d-2}\int_D\int_D
\frac{|y-z|^{\rho-d}dydz}{|x-y|^{d-2} |z-w|^{d-2}} \, +\,
|x-w|^{d-1}\int_D\int_D
\frac{|y-z|^{\rho-d}dydz} {|x-y|^{d-1}|z-w|^{d-2}}\\
&&~~~~+\, |x-w|^{d-1} \int_D\int_D
\frac{ |y-z|^{\rho-d}dydz} {|x-y|^{d-2}|z-w|^{d-1}}+\, |x-w|^{d}\int_D\int_D \frac{|y-z|^{\rho-d} dydz}
{|x-y|^{d-1}
|z-w|^{d-1}} \Big)
\end{eqnarray*}
for some constant $c=c(D,a)>0$.  Now combining the above with
\cite[Lemma 3.12]{GR},we easily get the conclusion of the proposition.
\qed

\begin{prop}\label{l:GR2d}
Suppose that $d\ge 1$ and that (A2) holds.
There exists a positive constant $C_{41}=C_{41}(D, a)$
 such that for every $(x,w) \in D \times D$
$$
\int_D\int_D  G^{a}_D(x,y)j^a_1(y-z)  G^{a}_D(z,w) dydz \,\le
\,C_{41}\,\frac{\delta_D(x)\delta_D(w)} {|x-w|^{d-\rho}}.
$$
\end{prop}

\pf
Recall that for every $d\ge 1$, $G^{a}_D(x,y) \le c_1 \frac{\delta_D(x)
\delta_D(y)}{|x-y|^d}$ by Theorem \ref{t-main-green}.
Therefore, by
following the arguments in the proof of \cite[Lemma 8]{R} we have
\begin{equation}\label{Rf1}\int_D
G^{a}_D(x,y)\frac{1}{|y-z|^{d-\rho}}\, dy\leq c_2
\frac{\delta_D(x)}{|x-z|^{d-\rho}}.
\end{equation}
Thus
\begin{eqnarray*}
\int_D G^{a}_D(x,y)j^a_1(y-z)dy\leq c_3
\frac{\delta_D(x)}{|x-z|^{d-\rho}}
\end{eqnarray*}
and so, since $G^{a}_D(z,w)=G^{a}_D(w,z)$, by \eqref{Rf1}
$$
\int_D\int_D  G^{a}_D(x,y)j^a_1(y-z)  G^{a}_D(z,w) dydz
\,\le \,c_3\,\delta_D(x) \int_D  G^{a}_D(z,w)
\frac{1}{|x-z|^{d-\rho}}   dz \le c_4\frac{\delta_D(x)\delta_D(w)}
{|x-w|^{d-\rho}}.
$$
\qed

\begin{lemma}\label{potlematuzup:this}
Suppose that $d\in \{1, 2\}$,
 $T>0$ and that (A1) holds. Then there exists a constant
$C_{42}=C_{42}(a,T)>0$ such that
$$
\sup_{0<t\le T}  \left(p^a(t,x,y)-e^{-2{\cal J}_at}p^Z(t,x,y)\right)
\le C_{42}
t^{1-d/2}.
$$

\end{lemma}
\pf
Recall that  $p^Z(t,y-x):= p^Z(t,x,y)$ and $p^a(t,y-x)=p^a(t,x,y)$.
Let $\wh Z$ be a pure jump L\'evy process with L\'evy density
$j^a_1$ in $\R^d$ independent of $Z$. Then $\wh Z$ is a compound
Poisson process with transition probability given by
$$
P^{\wh Z}(t, \cdot)=e^{-{\cal J}_at}\delta_0(\cdot)+ e^{-{\cal
J}_at}\sum^{\infty}_{n=1}\frac{t^n (j^a_1)^{*n}(\cdot)}{n!}.
$$
The process $Z+\wh Z$ has the same distribution as $X^a$. Thus the
distribution of $X^a_t$ is equal to the convolution of $p^Z(t,
\cdot)$ and $P^{\wh Z}(t, \cdot)$. Consequently, we have
$$
p^a(t, x)=p^Z(t, x)e^{-{\cal J}_at}+ e^{-{\cal
J}_at}\sum^{\infty}_{n=1} \frac{t^n p^Z(t,
\cdot)*(j^a_1)^{*n}(x)}{n!}.
$$
It follows from Lemma \ref{gestzabity1:this} and \eqref{e:CK} that
for $0<t \le T$
$$
p^Z(t,\cdot)*(j^a_1)^{*n} (x)\le  e^{{\cal J}_at}
p^a(t,\cdot)*(j^a_1)^{*n} (x)\le c_1 e^{{\cal J}_at}({\cal J}_a)^n
(t^{-d/\alpha} \wedge t^{-d/2})
\leq c_2 t^{-d/2} e^{{\cal J}_at}({\cal J}_a)^n
$$
for some positive constants $c_1, c_2$.
Thus it follows from Lemma
\ref{gestzabity1:this}
 and the above two displays that
for $0<t \le T$
\begin{eqnarray*}
p^a(t, x)- e^{-2{\cal J}_at} p^Z(t, x)&=& p^a(t, x)- e^{-{\cal
J}_at} p^Z(t, x) + e^{-{\cal J}_at}(1-e^{-{\cal J}_at})p^Z(t, x)\\
&\le &
c_2\sum^{\infty}_{n=1}
\frac
{t^{n-d/2}({\cal J}_a)^{n}}{n!} +
(1-e^{-{\cal J}_at})p^a(t, x).
\end{eqnarray*}
Since by \eqref{e:CK} $p^a(t, x)\leq
c_3 t^{-d/2}$
for $0<t \le T$, we reach the
conclusion of the lemma in view of the above display. \qed

We also need the following lemma.

\begin{lemma}\label{L:6.12}
Let $D$ be a $C^{1,1}$ open set with $C^{1,1}$ characteristics $(R,
\Lambda)$. Then there is a constant $C_{43}>0$ such that for all
$x\in D$ and $t>0$,
$$
\P_x(\tau^a_D>t) \,\leq \, C_{43}\, \frac{\delta_D(x)}{t+\sqrt{t}}.
$$
\end{lemma}

\pf When $t\geq 1$, the above inequality follows immediately from
Markov's inequality and  \eqref{e:g01}.  To establish the inequality
for the case of $0<t<1$, we will use a result from \cite{CKSV}.

We will only give the proof for the case $d\ge 2$. The proof in the
case $d=1$ is similar but simpler. Without loss of generality, we
can always assume that $R\leq 1$ and $\Lambda \geq 1$. By
definition, for every $Q\in \partial D$, there is a
$C^{1,1}$-function $\phi_Q: \bR^{d-1}\to \bR$ satisfying
$\phi_Q(0)=0$, $ \nabla\phi_Q (0)=(0, \dots, 0)$, $\| \nabla \phi_Q
\|_\infty \leq \Lambda$, $| \nabla \phi_Q (x)-\nabla \phi_Q (z)|
\leq \Lambda |x-z|$,
and an orthonormal coordinate system $CS_Q$: $y=(\wt y, \, y_d)$
such that $B(Q, R )\cap U=\{ y=(\wt y, \, y_d)\in B(0, R) \mbox{ in
} CS_Q: y_d > \phi (\wt y) \}$. Define
$$
\rho_Q (x) := x_d -  \phi_Q (\wt x),
$$
where $(\wt x, x_d)$ are the coordinates of $x$ in $CS_Q$. Note that
for every $Q \in \partial U$ and $ x \in B(Q, R)\cap U$, we have
$
(1+\Lambda^2)^{-1/2} \rho_Q (x) \le \delta_U(x) \le
\rho_Q(x).$
We define for $r_1, r_2>0$
$$
D_Q( r_1, r_2) :=\left\{ y\in U: r_1 >\rho_Q(y) >0,\, |\wt y | < r_2
\right\}.
$$
Note that for $b>0$,
\begin{eqnarray*}
\P_x(\tau^{b}_{D}>1)& \leq&
\P_x \left(\tau^{b}_{ D_Q( \delta_0 , r_0)}>1  \right)+\P_{x}
\left(X^{b}_{\tau^{b}_{ D_Q( \delta_0,  r_0)}}
\in  D \hbox{ and } \tau^{b}_{ D_Q( \delta_0,  r_0)} \leq 1 \right)\\
&\le & E_x\left[\tau^b_{  D_Q( \delta_0,  r_0)}\right]+
\P_{x}\left(X^{b}_{\tau^{b}_{ D_Q( \delta_0,  r_0)}} \in  D\right).
\end{eqnarray*}
Thus  by
\cite[Lemma 3.5]{CKSV}, there is a constant
$c_1=c_1(R, \Lambda, a)$ so that for every $b\in (0, a]$
\begin{equation}\label{e:6.8}
\P_x(\tau^{b}_{D}>1) \le c_1 \delta_D (x) \qquad \hbox{for every }  x\in D.
\end{equation}
Note that for $0<\lambda\leq 1$, $\lambda^{-1}D$ is a $C^{1,1}$ open
set with $C^{1, 1}$ characteristics $(R, \Lambda)$. Hence by the
scaling property of $X^a$
in \eqref{e:scaling}, we have from \eqref{e:6.8} that for $t\in (0,
1]$,
\begin{eqnarray*}
\P_x (\tau^a_D>t)
\ =\
 \ \P_{t^{-1/2}x} \left( \tau^{a t^{(2-\alpha)/(2\alpha)}}_{t^{-1/2}D}>1\right)
\ \leq \ c_1 \ \delta_{t^{-1/2}D}(t^{-1/2}x) \ =  \ c_1 \ \frac{\delta_D(x)}{\sqrt{t}}.
\end{eqnarray*}
This completes the proof of the lemma. \qed

\begin{thm}\label{t:gee}
Suppose that (A1)--(A3) hold. There exists $C_{44}=C_{44}(D, a)
>0$ such that
\begin{equation}\label{e:ge}
C_{44}^{-1} \,G^Z_D(x,w) \,\le\, G^a_D(x,w) \,\le\, C_{44}\, G^Z_D(x,w), \quad
(x,w) \in D \times D.
\end{equation}
\end{thm}

\pf By \eqref{e:G_ub} and Lemma \ref{G:g2}, we only need to show the
second inequality in (\ref{e:ge}) for $|x-y|^2<\eps$, where
$\eps\in (0, 1)$ is a constant to be chosen later.  We consider the
cases
$d\ge3$ and $d\le 2$ separately.

\medskip \noindent
(a) $d\ge3$:
Applying Lemma \ref{l:g5}  and then Proposition \ref{l:GR2}, we get
$$
G^{a}_D(x,y) \, \le \,  G^Z_D(x,y) \,+\, c_1 \, G^{a}_D(x,y)
|x-y|^{d-2}.
$$
Choose $\eps>0$ small so that
$$
c_1 \, G^{a}_D(x,y) |x-y|^{d-2} \, \le \, \frac12  G^{a}_D(x,y)
\quad \mbox{ if } |x-y| < \eps^{1/2}.
$$
Thus
$$
 G^{a}_D(x,y) \, \le \, 2  G^Z_D(x,y)  \quad \mbox{ if } |x-y| < \eps^{1/2}.
$$

\medskip \noindent

(b) $d\le 2$: We first note that, since $p^{a}_D(t,x,y) \le c_2
(t^{-d/\alpha} \wedge t^{-d/2})$ by \eqref{e:CK}, using the
semigroup property,
\begin{eqnarray*}
p^{a}_D(t, x, y) &=& \int_D p^{a}_D(t/3, x, z) \int_D p^{a}_D(t/3,
z, w) p^{a}_D(t/3, w, y) dw dz\\
&\le& c_2
(t^{-d/\alpha} \wedge t^{-d/2})
 \int_D p^{a}_D(t/3, x, z)
dz\int_D p^{a}_D(t/3, w, y) dw\\
&=& c_2 (t^{-d/\alpha} \wedge t^{-d/2})\, \P_x ( \tau^{a}_D > t/3)
\P_y (  \tau^{a}_D > t/3).
\end{eqnarray*}
By Lemma \ref{L:6.12}, we get
\begin{eqnarray*}
p^{a}_D(t, x, y) &\le& c_3 \left(t^{-d/\alpha} \wedge t^{-d/2}\right)
 \, \left(
\frac{\delta_D(x)}{\sqrt{t}} \wedge 1\right)
 \left( \frac{\delta_D(y)} {\sqrt{t}} \wedge 1\right) \\
&\le& c_4 \left(t^{-d/\alpha-1}\wedge  t^{-d/2-1}\right)
 \delta_D(x) \delta_D(y) .
\end{eqnarray*}
Consequently,
\begin{equation}\label{e:qwe}
\int_{t_0}^\infty p^{a}_D(t, x, y)\, dt \le c_5t_0^{-d/2}\,
\delta_D(x)\delta_D(y)
\qquad \hbox{for every } t_0>0 \hbox{ and } x, y\in D.
\end{equation}
Since
$$
p^{a}_D(t,x,y) = p^{a}(t,x,y)-\E_x\left[p^{a}(t-\tau^{a}_D,
X^{a}_{\tau^{a}_D},y):t\geq \tau^{a}_D \right]
$$
and
$$
p^Z_D(t,x,y) = p^Z(t,x,y)-\E_x\left[p^Z(t-\tau^Z_D,Z_{\tau^Z_D},y): t\geq
\tau^Z_D \right],
$$
it follows from \eqref{e:new} that
\begin{eqnarray*}
p^{a}_D(t,x,y)&=&
p^{a}(t,x,y)-\E_x\left[p^{a}(t-\tau^{a}_D,
X^{a}_{\tau^{a}_D},y):t\geq \tau^{a}_D \right]
\\
&\le&
p^{a}(t,x,y)-e^{-2{\cal
J}_at} \E_x\left[p^Z(t-\tau^Z_D,Z_{\tau^Z_D},y): t\geq
\tau^Z_D \right]\\
&=&
p^{a}(t,x,y)-e^{-2{\cal
J}_at} (p^Z(t,x,y)-p_D^Z(t,x,y))\\
&\leq& p^{a}(t,x,y)+p^Z_D(t,x,y)-e^{-2{\cal
J}_at}p^Z(t,x,y).
\end{eqnarray*}
So integrating over $[0,t_0]$ with $t_0= (\delta_D(x)
\delta_D(y))^{1/2}$, which is bounded by ${\rm diam}(D)$, we have by
Lemma \ref{potlematuzup:this} and \eqref{e:qwe} that
\begin{eqnarray}
G^{a}_D(x,y)&=& \int_{0}^{t_0}p^{a}_D(t,x,y)dt+ \int^{\infty}_{t_0}
p^{a}_D(t,x,y)dt\notag\\
&\leq& G^Z_D(x,y)+ \int^{t_0}_0(p^{a}(t,x,y)-e^{-2{\cal
J}_at}p^Z(t,x,y))dt+c_5
t^{-d/2}_{0} \delta_D(x)\delta_D(y)\notag\\
&\leq& G^Z_D(x,y) + c_6 t_0^{2-d/2} +c_5
t^{-d/2}_{0} \delta_D(x)\delta_D(y)\notag\\
&\leq& G^Z_D(x,y) + c_7(\delta_D(x)\delta_D(y))^{1-d/4}.
\label{Gfala3:this2}
\end{eqnarray}
Since $G^{a}_D(x,y)\ge c_8(\delta_D(x)\delta_D(y))^{1-d/2}$ for $|x-y|^2\leq\delta_D(x)\delta_D(y)$
by Theorem \ref{t-main-green},
 we have from(\ref{Gfala3:this2}) that
\begin{equation}\label{deltaEstimate}
G^{a}_D(x,y)\leq G^Z_D(x,y)+ c_9 (\delta_D(x)\delta_D(y))^{d/4}
G^{a}_D(x,y) \quad \text{ if } |x-y|^2\leq\delta_D(x)\delta_D(y).
\end{equation}
On the other hand, applying Lemma \ref{l:g5}  and then Proposition
\ref{l:GR2d}, we get
$$
G^{a}_D(x,y) \, \le \,  G^Z_D(x,y) \,+\, c_{10} \, \frac{\delta_D(x)
\delta_D(y)} {|x-y|^{d-\rho}}.
$$
Since $c_{11}  \frac{\delta_D(x)\delta_D(y)} {|x-y|^{d}} \le
G^{a}_D(x,y)$ for  $|x-y|^2\geq\delta_D(x)\delta_D(y)$  by Theorem
\ref{t-main-green},
we have
\begin{equation}\label{deltaEstimate2}
G^{a}_D(x,y) \, \le \,  G^Z_D(x,y) \,+\, c_{12} \,
|x-y|^{\rho}G^{a}_D(x,y).
\end{equation}
Now using \eqref{deltaEstimate}-\eqref{deltaEstimate2}, we can
choose $\eps\in (0, 1)$ small so that
$$
c_{12} \, G^{a}_D(x,y) |x-y|^{\rho} \, \le \, \frac12 G^{a}_D(x,y)
\quad \mbox{ if } \delta_D(x)\delta_D(y) \le |x-y|^2 < \eps
$$
and
$$
c_9 (\delta_D(x)\delta_D(y))^{d/4}G^{a}_D(x,y) \le \, \frac12
G^{a}_D(x,y) \quad \mbox{ if } |x-y|^2 \le \delta_D(x)\delta_D(y) <
\eps.
$$
Thus in these cases,
$
 G^{a}_D(x,w) \, \le \, 2  G^Z_D(x,w).$

For the remaining case $\delta_D(x)\delta_D(y)\ge  \eps$, we use
\eqref{Gfala3:this2}, Lemma \ref{momentyogr:this} and \eqref{lbG} to
get that
$$
G^{a}_D(x,w) \, \le (1 + c_{13}(\delta_D(x)\delta_D(y))^{-d/4})
G^Z_D(x,y) \le (1 + c_{13}\eps^{-d/4}) G^Z_D(x,y).
$$
The proof of the theorem is now complete. \qed

We now show that the theorem above covers the case of the
independent sum of a Brownian motion and a relativistic stable
process, and the case of the independent sum of a Brownian motion
and a truncated stable process.

For any $m\ge 0$, a relativistic $\alpha$-stable process $Y^m$
in $\bR^d$ with mass $m$ is a L\'evy process with characteristic
function given by
$$
\E_x \left[e^{i \xi \cdot (Y^m_t-Y^m_0) } \right] = \exp\left(-t
\left( \left(|\xi|^2+ m^{2/\alpha} \right)^{\alpha/2}-m\right)
\right), \qquad \xi \in \bR^d.
$$
Suppose that $Y^m$ is independent of the Brownian motion $X^0$.
We define $Z^m$ by $Z^m_t:=X^0_t+ Y^m_t$.  We will call the process
$Z^m$ the independent sum of the Brownian motion $X^0$ and the
relativistic $\alpha$-stable process $Y^m$ with mass $m$. The
L\'{e}vy measure of $Z^m$ has a density
$$
J^{Z^m}(x)=  {\cal A} (d, \, -\alpha) |x|^{-d-\alpha} \psi
(m^{1/\alpha}|x|)
$$
where
$$
\psi (r):= 2^{-(d+\alpha)} \, \Gamma \left(
\frac{d+\alpha}{2} \right)^{-1}\, \int_0^\infty s^{\frac{d+\alpha}{
2}-1} e^{-\frac{s}{ 4} -\frac{r^2}{ s} } \, ds,
$$
which is a decreasing smooth function of $r^2$. (see \cite[pp.
276-277]{CS4} for  details). Thus
$$
0 \le j^1(x)-j^{Z^m}(x) \le c |x|^{2-\alpha-d}.
$$
Moreover, the conditions
(A1), (A2) and (A3)(a)
can be checked easily.
Therefore as a corollary of Theorem \ref{t:gee}, we have the
following.
\begin{cor}\label{c:GM1}
There exists a constant $C_{45}=C_{45}(D, \alpha)>0$ such that
$$
C_{45}^{-1}G^1_D(x, y)\le G^{Z^m}_D(x, y)\le C_{45} G^1_D(x, y),
\quad x, y \in D,
$$
where $G^{Z^m}_D(x, y)$ is the Green function of
$Z^m$ in $D$.

\end{cor}

By a $\lambda$-truncated symmetric $\alpha$-stable process in $\R^d$
we mean a pure jump symmetric L\'evy process
$\wh Y^\lambda=(\wh Y^\lambda_t, t \ge 0, \P_x, x\in \R^d )$
in $\R^d$  with the  L\'evy density ${\cal A}(d, -\alpha)
|x|^{-d-\alpha}\, \1_{\{|x|< \lambda\}}$.  Note that the
L\'evy exponent $\psi^\lambda$ of $\wh Y^\lambda$, defined by
$$
\E_x\left[e^{i\xi\cdot(\wh Y^\lambda_t-\wh Y^\lambda_0)}\right] =
e^{-t\psi^\lambda(\xi)} \quad \quad \mbox{ for every } x\in \R^d
\mbox{ and } \xi\in \R^d,
$$
is given  by
\begin{equation}\label{e:psi}
\psi^\lambda(\xi)={\cal A}(d, -\alpha) \int_{\{|y|<
\lambda\}}\frac{1-\cos(\xi\cdot y)}{|y|^{d+\alpha}}dy.
\end{equation}

Suppose that $\wh Y^{\lambda}$ is a $\lambda$-truncated
symmetric $\alpha$-stable process in $\R^d$ which is independent of
the Brownian motion $X^0$. We define
$
\wh X_t^{\lambda}:=X^0_t+ \wh Y^{\lambda}_t$ for $t\geq 0 .
$
Then $\wh X^{ \lambda}$ has the same distribution as the
L\'evy process obtained from $X^1$ by removing jumps of size larger
than $\lambda$.

Unlike the symmetric stable process
$Y$, the process $\wh Y^{\lambda}$ can only
make jumps of size less than
$\lambda.$ In order to guarantee
the strict positivity of the transition density
$p^{\wh X^{\lambda}}_D(t,x,y)$ for $\wh X^{\lambda, D}$,
we need to
 impose
 the following assumption on $D$.

\begin{defn}\label{d:rc}
We say that an open set $D$ in $\R^d$ is $\lambda$-roughly connected if for
every $x, y \in D$, there exist finitely many distinct connected
components $U_{1}, \cdots, U_{m}$ of  $D$ such that $ x \in U_{1}$,
$y \in U_{m}$ and dist$(U_{k}, U_{{k+1}}) <\lambda$ for $1 \le k \le
m-1$.
\end{defn}

The following result is proved in \cite{KS6}.

\begin{prop}
[{\cite[Proposition 4.4]{KS6}}]
\label{p:spdp}
For any bounded $\lambda$-roughly connected open set $D$ in $\R^d$,
the transition density
$p^{\wh X^{\lambda}}_D(t,x,y)$ for $\wh X^{\lambda, D}$
is strictly positive in $(0, \infty)\times D \times
D$.
\end{prop}
The other conditions (A1), (A2) and (A3)(b) can be checked easily.
Therefore as a corollary of Theorem \ref{t:gee},
 we have the following.
\begin{cor}\label{c:GM2}
Suppose $D$ is a bounded $\lambda$-roughly connected $C^{1, 1}$ open
set in $\bR^d$,
$d\ge 1$.
There exists $C_{46}=C_{46}(D, \alpha)>0$
such that
$$
C_{46}^{-1}G^1_D(x, y)\le
G^{\wh X^{\lambda}}_D(x, y)
\le C_{46} G^1_D(x, y), \quad x, y \in D,
$$
where
$G^{\wh X^{\lambda}}_D(x, y)$ is the Green function of
$\wh X^{\lambda}$ in $D$.

\end{cor}

\vskip 0.3truein

{\bf Zhen-Qing Chen}

Department of Mathematics, University of Washington, Seattle,
WA 98195, USA

E-mail: \texttt{zchen@math.washington.edu}

\bigskip

{\bf Panki Kim}

Department of Mathematical Sciences and Research Institute of Mathematics,
Seoul National University,
San56-1 Shinrim-dong Kwanak-gu,
Seoul 151-747, Republic of Korea

E-mail: \texttt{pkim@snu.ac.kr}

\bigskip

{\bf Renming Song}

Department of Mathematics, University of Illinois, Urbana, IL 61801, USA

E-mail: \texttt{rsong@math.uiuc.edu}

\bigskip

{\bf
Zoran Vondra\v{c}ek}

Department of Mathematics,
University of Zagreb,
Bijeni\v{c}ka c.~30,
Zagreb, Croatia

Email: \texttt{vondra@math.hr}

\end{document}